\documentclass[12pt,reqno]{amsart}
\usepackage{amsthm,amsfonts,amssymb,amscd}
\usepackage[all]{xy}

\newtheorem{theorem}[subsection]{Theorem}

\newtheorem{proposition}[subsection]{Proposition}

\newtheorem{lemma}[subsection]{Lemma}
\newtheorem{corollary}[subsection]{Corollary}

\theoremstyle{definition}
\newtheorem{definition}[subsection]{Definition}

\newtheorem{proposition-definition}[subsection]{Proposition-Definition}

\theoremstyle{remark}
\newtheorem{remark}[subsection]{Remark}

\newcommand{\actson}{\begin{picture}(7,4)\thinlines\linethickness{0.02mm}
\qbezier(2,-0.2)(8,1.3)(2,2.8)
\put(2,2.8){\vector(-2,1){0.5}}\end{picture}
}
\newcommand{\ccd}{\begin{picture}(2,2)
\put(0.5,1){\circle*{0.7}}\end{picture}}
\newcommand{\gseven}{\Gamma_{12}^{7}}

\newcommand{\dual}{{\scriptscriptstyle \vee}}
\newcommand{\sten}{{\Sigma_{12}^{10}}}

\newcommand{\coker}{\operatorname{coker}\nolimits}

\newcommand{\elm}{\operatorname{elm}\nolimits}
\newcommand{\END}{\operatorname{End}\nolimits}

\newcommand\Hilb{{\operatorname{Hilb}\nolimits}}

\newcommand{\length}{\operatorname{length}\nolimits}

\newcommand{\mult}{\operatorname{mult}\nolimits}

\newcommand{\pf}{\operatorname{Pf}\nolimits}
\newcommand{\Pic}{\operatorname{Pic}\nolimits}
\newcommand{\pr}{\operatorname{pr}\nolimits}
\newcommand{\red}{\operatorname{red}}

\newcommand{\rk}{\operatorname{rk}\nolimits}

\newcommand{\Sing}{\operatorname{Sing}}

\newcommand{\spin}{\operatorname{Spin}(8)}
\newcommand{\Spin}{\operatorname{Spin}(10)}

\newcommand\Sym{{\operatorname{Sym}\nolimits}}

\newcommand{\CC}{{\mathbb C}}
\newcommand{\RR}{{\mathbb R}}
\newcommand{\ZZ}{{\mathbb Z}}

\newcommand{\PP}{{\mathbb P}}

\newcommand{\OOO}{{\mathcal O}}

\newcommand{\III}{{\mathcal I}}

\newcommand{\EEE}{{\mathcal E}}

\newcommand{\FFF}{{\mathcal F}}

\newcommand{\NNN}{{\mathcal N}}

\newcommand{\QQQ}{{\mathcal Q}}

\newcommand{\TTT}{{\mathcal T}}

\newcommand{\UUU}{{\mathcal U}}
\newcommand{\VVV}{{\mathcal V}}

\newcommand{\EEEE}{\boldsymbol{\mathcal E}}

\newcommand{\cal}[1]{\mathcal{#1}}
\renewcommand{\bar}[1]{\overline{#1}}
\newcommand{\equi}{\: \Longleftrightarrow\: }

\newcommand\alp{\alpha}

\renewcommand\phi{\varphi}

\newcommand\Si{\Sigma}

\newcommand{\into}{\hookrightarrow}

\newcommand\lra{{\longrightarrow}}

\renewcommand\empty{\varnothing}

\author{A. Iliev}

\address{\scriptsize Atanas Iliev: 
Institute of Mathematics,
Bulgarian Academy of Sciences,
Acad. G. Bonchev Str., 8,\ 
1113 Sofia, Bulgaria}
\email{ailiev@math.bas.bg}

\author{D. Markushevich}\thanks{D.M.: Partially supported by the grant
INTAS-OPEN-2000-269}

\address{\scriptsize Dimitri Markushevich:  
Math\'ematiques - b\^{a}t.M2, Universit\'e Lille 1, 
F-59655 Villeneuve d'Ascq Cedex, France}
\email{markushe@agat.univ-lille1.fr}


\title{Elliptic curves and rank-2 vector bundles\\
on the prime Fano threefold of genus 7}

\begin{document}

\begin{abstract}
According to Mukai, any prime
Fano threefold $X$ of genus 7 is a linear section of the
spinor tenfold in the projectivized half-spinor space
of Spin(10). It is proven that the moduli space of stable
rank-2 vector bundles with Chern classes
$c_1=1,c_2=5$ on a generic X is isomorphic to the
curve of genus 7 obtained by taking an orthogonal
linear section of the spinor tenfold. This is an inverse
of Mukai's result on the isomorphism of a non-abelian Brill--Noether
locus on a curve of genus 7 to a Fano threefold
of genus 7. An explicit geometric construction
of both isomorphisms and a similar result for K3 surfaces
of genus 7 are given.
\end{abstract}

\maketitle

\setcounter{section}{-1}

\section{Introduction}

The study of moduli spaces of stable vector bundles
on Fano threefodls of indices 1 and 2 is quite a
recent topic. The {\em index} of a Fano threefold $X$ is the maximal
integer $\nu$ dividing $K_X$ in the Picard group of $X$.
The known results sofar include
the desciption of one component of moduli of
rank-2 vector bundles on each one of the following four Fano threefolds:
the cubic \cite{MT}, \cite{IM-1}, \cite{Druel}, \cite{B},
the quartic \cite{IM-2}, the prime Fano threefold
of genus 9 \cite{IR}
and the double solid of index 2 \cite{Ti}.
It turns out, that the flavour of the results one can obtain
depends strongly on the index.
In the index-2 case, the answers are given in terms
of the Abel--Jacobi map of the moduli of vector bundles
into the intermediate Jacobian $J(X)$, defined by the
second Chern class $c_2$, and the techniques
originate to Clemens, Griffiths, and Welters.
For the cubic $X_3$, the moduli space $M_{X_3}(r;c_1,c_2)$
with invariants $r=2, c_1=0,c_2=5$ is of dimension 5
and is identified with an open subset of the intermediate Jacobian.
For the double solid $Y_2$ of index 2, Tikhomirov
found a 9-dimensional component of $M_{Y_2}(2;0,3)$ whose Abel--Jacobi
map is quasi-finite onto an open subset of the theta-divisor
of $J(Y_2)$, and it is conjectured and almost proven that this map
is in fact birational (in preparation).

In the index-1 case, the Abel--Jacobi map does not bring
much new information about the moduli spaces investigated
up to now. For the quartic threefold $X_4$, we proved in
\cite{IM-2} that $M_{X_4}(2;1,6)$ has a component of dimension 7
with a 7-dimensional Abel--Jacobi image in the 30-dimensional
intermediate Jacobian $J(X_4)$. One can conclude from here about the
geometry of this component only that its Kodaira dimension
is positive.

In the present paper, we consider one more index-1 case:
we determine the moduli space
$M_X=M_X(2;1,5)$ for a generic prime Fano threefold $X=X_{12}$ of genus 7.

Following the classical terminology, we call
the Fano threefolds $X_{2g-2}$ of index 1
and degree $2g-2$ with Picard number 1
{\em prime Fano threefolds of genus} $g$.
They have been classified, up to deformations, by
Iskovskikh \cite{I1}, \cite{IP}. There is only one
moduli family of threefolds $X_{2g-2}$ for every
$g=2,\ldots, 12, g\neq 11$ (see Table 12.2 in \cite{IP}).
Mukai \cite{Mu-1} proved that
$X_{2g-2}$ is a linear section of
some projective homogeneous space $\Sigma_{2g-2}$ for $7\leq g\leq 10$. 
In the case $g=7$, $\Sigma =\Sigma_{12}$ is
the spinor tenfold in $\PP^{15}$.
It is self-dual, that is, the dual variety $\check{\Sigma}\subset
\check{\PP}^{15}$, formed by the hyperplanes in $\PP^{15}$
tangent to $\Sigma$, is isomorphic to $\Sigma$ via some
projectively linear map identifying $\PP^{15}$ with its
dual $\check{\PP}^{15}$. Thus, to a linear section
$X=\PP^{7+k}\cap\Sigma$ of $\Sigma$ of dimension $2-k$ we can
associate the orthogonal linear section
$\check{X}:=\check{\PP}^{7-k}\cap\check{\Sigma}$ of dimension $2+k$,
where $\check{\PP}^{7-k}=(\PP^{7+k})^{\perp}\subset \check{\PP}^{15}$.
For $k=1$, we obtain a curve linear section $\Gamma =\check{X}$,
which is a canonical curve of genus 7. Our main result is the
following statement.

\begin{theorem}\label{main}
Let $X=X_{12}$ be a generic prime Fano threefold of genus $7$.
Then $M_X$ is isomorphic to the curve $\Gamma =\check{X}$.
\end{theorem}

We prove also similar statements in the cases $k=0$, where
$X,\check{X}$ are generic K3 surfaces of degree 12, and
$k=-1$, where $X$ is a curve and $\check{X}$ is a threefold.
In the latter case, one should take the non-abelian Brill-Noether
locus of rank-2 vector bundles on $X$ with canonical determinant
and 5 linearly independent global sections on the role of $M_X$.
For $k=0$, Mukai \cite{Mu-0}, \cite{Mu-4} proved that $M_X$
is another K3 surface
of degree 12 ($M_X$ represents the so-called Fourier--Mukai transform of $X$;
see \cite{HL}). We make this statement more precise
by identifying $M_X$ with the orthogonal K3 surface
$\check{X}$ via an explicit map $\rho_X$ having a beautiful geometric
construction.

For $k=-1$, Mukai \cite{Mu-5} proved that $M_X$
is a Fano threefold of degree 12. Again, we show
that this Fano threefold is isomorphic to the orthogonal
linear section of the spinor tenfold, and our Main Theorem
represents the inverse of this result.

Iliev--Ranestad \cite{IR}
obtained similar results for the 1-, 2- and 3-dimensional
linear sections of the symplectic Grassmannian $\Sigma_{16} 
\subset \PP^{13}$, but in their case the dual of $\Sigma_{16}$
is a singular quartic hypersurface in $\check{\PP}^{13}$,
so the moduli spaces
(or the non-abelian Brill--Noether locus in dimension 1)
that they consider are isomorphic to linear sections
of this quartic hypersurface.

Our construction of the map $\rho_X:\check{X}\lra M_X$
is very simple: for any $w\in\check{X}$ the corresponding
hyperplane $\PP^{14}_w$ in $\PP^{15}$ is tangent to $\Sigma$
along a projective space $\PP^4$, and the linear projection $\pi_w$
of $X\cap\PP^{14}_w$ into $\PP^9$ with center $\PP^4$
has its image inside the Grassmannian $G(2,5)\subset\PP^9$.
It turns out that the pullback of the universal rank-2 bundle
from $G(2,5)$ to $X$ is stable
and its class belongs to $M_X$. This defines the
image of $w$ in $M_X$. It is not so obvious that
the thus defined map $\rho_X$ is nontrivial. In Proposition \ref{image-open}
we prove that its image is an irreducible component of $M_X$.
This enables us to conclude the proof of the fact that $\rho_X$
is an isomorphism in the cases, where the irreducibility
and smoothness of $M_X$ are already known from the work of Mukai:
$k=-1$ (Proposition \ref{irred-curve}) and $k=0$ (Proposition \ref{irredK3}).
For the case $k=1$, we prove in Proposition \ref{MXglobgen}
that every vector bundle $\EEE\in M_X$ is globally generated
and is obtained by Serre's construction from a normal
elliptic quintic contained in $X$. The irreducibility
of $M_X$, equivalent to that of the family of elliptic quintics
in $X$, is reduced to the known irreducibility
in the K3 case.  The smoothness of $M_X$ is proved separately in using
Takeuchi--Iskovskikh--Prokhorov birational maps
$\Phi_p:X\dasharrow Y=Y_5$ and $\Psi_q:X\dasharrow Q$, where
$p\in X$ is a point, $q\subset X$ a conic, $Y_5$
the Del Pezzo threefold of degree 5 and  $Q$
the three-dimensional quadric hypersurface.

Takeuchi \cite{Tak} has undertaken a systematic study
of birational transformations of Fano varieties
that can be obtained by a blow up with
center in a point $p$, a line $\ell$ or a conic $q$ followed by a flop
and a contraction of one divisor. Iskovskikh--Prokhorov \cite{IP}
have extended Takeuchi's list,
in particular, they found the two birational transformations
for $X_{12}$ mentioned above.
The techniques of proofs are those of Mori theory,
based on the observation that the Mori cone of $X_{12}$
blown up at a point, line or conic is an angle
in $\RR^2$, hence there are exactly two extremal rays
to contract, the first one giving the initial 3-fold,
the second one defining the wanted birational
map. But before one can contract the second extremal
ray, one has to make a flop.
We describe in detail the structure of $\Phi_p$, $\Psi_q$
(Theorems \ref{bir1} and \ref{bir2}). The last contraction
in both cases blows down one divisor onto a curve
of genus 7. Thus, we have 3 curves of genus 7 associated to $X$:
the orthogonal linear section $\Gamma$, and $\Gamma ',\Gamma ''$ coming
from the birational maps. We prove that the three curves are
isomorphic.
We also identify the flopping curves for $\Phi_p$: they are
the 24 conics passing through $p$, and their images are
the 24 bisecants of $\Gamma '$.

The ubiquity of the
maps $\Phi_p$, $\Psi_q$ is in that they provide a stock of well-controlled
degenerate elliptic quintics: the ones with a node at $p$ are just the
proper transforms of the unisecant lines of $\Gamma'$ in $Y$
and the reducible ones having $q$ as one of components are nothing else but
the proper transforms of the exceptional curves contracted
by $\Psi_q$ into points of $\Gamma ''$, that is,
they are parametrized by $\Gamma ''$.
The smoothness of $M_X$ follows from the existence, among the zero loci
of sections of any vector bundle $\EEE\in M_X$, of a nodal quintic
with a node at $p$ such that the normal bundle of the corresponding
unisecant of $\Gamma '$ is $\OOO\oplus\OOO$ (see the proofs
of Proposition \ref{smoothness} and Lemma \ref{ncx-nodal}).
The family of lines on $Y$ is well known (see for example \cite{Ili}, \cite{FN}).
In particular, $Y$ contains a rational curve $C_6^0$ which is
a locus of points $z$ such that there is a unique line in $Y$ passing through
$z$,
and the normal bundle of this line is $\OOO (1)\oplus\OOO (-1)$.
Hence our proof of smoothness does not work in the case
when $\Gamma '$ meets $C_6^0$. We prove that in this case a generic
deformation $\Gamma_t'$ of $\Gamma '$ does not meet $C_6^0$ (Lemma \ref{DP2})
and that $\Gamma_t'$ corresponds to some birational map
$\Phi_{p_t}:X_t\dasharrow Y$ of the same type. This explains the fact, why we
state our Main Theorem only for generic~$X$.
We conjecture that the conclusion
of the Theorem
is true for any smooth 3-dimensional linear section
$\Sigma\cap \PP^8$.

In Section \ref{tenfold}, we give a definition of the spinor tenfold $\Sigma$,
represent it as one of the two components of the family
of maximal linear subspaces of an 8-dimensional quadric,
and introduce the notion of the pure
spinor associated to a point of $\Sigma$.

In Section \ref{linsect}, we study some properties of linear sections $X$
of $\Sigma$, in particular, the projections $\pi_w$
to the Grassmannian $G(2,5)$ in $\PP^9$, defined by the points $w\in\check{X}$,
and prove that any linear embedding of $X$
into $G(2,5)$, under some additional restrictions,
is always given by such a projection (the case $\dim X=1$
is postponed until Section \ref{map-rho}).

In Section \ref{ellquint}, we list standard facts about
the moduli space $M_X$ ($\dim X =3$) and the Hilbert scheme
of elliptic quintics on $X$; we show that any $\EEE\in M_X$
is obtained by Serre's construction from a "quasi-elliptic"
quintic and that the fibers of Serre's construction
over $M_X$ are projective spaces $\PP^4$.

In Section \ref{map-rho}, we define the map $\rho_X:\check{X}\lra M_X$
in all the three cases $\dim X = 1,2,3$, and prove that its image
is a component $M_X^0$ of $M_X$. We prove that $\rho_X$
is an isomorphism for $\dim X=1,2$. For $\dim X =3$ we obtain the following
precision of the result of the previous section:
any $\EEE\in M_X$ is globally generated and
is obtained by Serre's construction from a {\em smooth} elliptic
quintic.

In the Sections \ref{bir-maps} and \ref{proof-MT}, $\dim X=3$.
In Section \ref{bir-maps}, we provide some basic properties of the families
of lines and conics on $X$, in particular, prove the irreducibility
of the family of conics, and describe the structure of the two
Takeuchi--Iskovskikh--Prokhorov birational maps. We show that
the vector bundles, constructed from the stock of quasi-elliptic quintics
generated by these maps, belong to $M_X^0$, and deduce the isomorphism
of the three curves $\Gamma ,\Gamma ',\Gamma ''$.

Section \ref{proof-MT} is devoted to the proof of Theorem \ref{main}.

\bigskip

{\sc Acknowledgements.} The second author thanks V.A.Iskovskikh,
who communicated to him Moishezon's Lemma \ref{moishe},
and Yu.~Prokhorov for discussions.

\section{Spinor tenfold}\label{tenfold}

The spinor tenfold 
$\sten$ is a homogeneous space of the complex spin group $\Spin$,
or equivalently, that of SO$(10)=\Spin /\{\pm 1\}$.
It can be defined as the unique closed orbit of $\Spin$
in the projectivized half-spinor representation $\Spin :\PP^{15}\actson$.
We will remind an explicit description of $\sten$ and
some of its properties, following essentially
\cite{C}, \cite{Mu-2}, \cite{RS}.

Let $Alt_5({\mathbb C}) \cong \mathbb C^{10}$ be the space of skew-symmetric
complex $5 \times 5$ matrices.
For $\hat{A} \in Alt_5({\mathbb C})$
denote by $\vec{\pf }(\hat{A}) \in \mathbb C^5$ the 5-vector with coordinates
$\vec{\pf }(\hat{A})_i = (-1)^i\pf _i(\hat{A})$, i = 1,\ldots ,5,
where $\pf _i$ are the codimension 1 Pfaffians of an odd-dimensional
skew-symmetric matrix.\bigskip

\begin{definition}\label{defsigma}
The 
Spinor $10$-fold $\Sigma = {\Sigma}^{10}_{12} \subset {\mathbb P}^{15}$ 
is the closure of the image $j(Alt_5({\mathbb C}))$ 
under the embedding 

\begin{equation}\label{alt5}
j:\CC^{10} \cong Alt_5(\mathbb C) \rightarrow \mathbb P^{15},\ 
\hat{A} \mapsto (1:\hat{A}:\vec{\pf }(\hat{A}))\end{equation} 
\end{definition}

We will write homogeneous coordinates in ${\mathbb P}^{15}$ 
in the form $(u:\hat{X}:\vec{y})$, where $u \in {\mathbb C}$, 
$\hat{X} = (x_{ij}) \in Alt_5({\mathbb C})$ 
and $\vec{y} = (y_1,...,y_5) \in {\mathbb C}^5$. 

The map $j$
parameterizes the points of the 
open subset $j(Alt_5({\mathbb C})) = \Sigma \cap (u \not= 0)$
and 

\begin{equation}\label{equa-sten}
(u:\hat{X}:\vec{y}) \in \Sigma \ \Longleftrightarrow \ 
              u\vec{y} = \vec{\pf }(X) \ \& \ \hat{X}\vec{y} = 0;
\end{equation}

\noindent 
Writing down the components of the above matrix equations,
we obtain the defining equations of $\Sigma$, or the generators of
the homogeneous ideal of $\Sigma \subset {\mathbb P}^{15}$:

\bigskip
\begin{figure}[!h]

\ \ \ \ \ \ \ \ \ \ \ \ \ \ \ \ \ \ \ \ $q_1^+ = uy_1 + x_{23}x_{45} -
x_{24}x_{35} + x_{34}x_{25}$

\ \ \ \ \ \ \ \ \ \ \ \ \ \ \ \ \ \ \ \ $q_2^+ = uy_2 - x_{13}x_{45} +
x_{14}x_{35} - x_{34}x_{15}$

\ \ \ \ \ \ \ \ \ \ \ \ \ \ \ \ \ \ \ \ $q_3^+ = uy_3 + x_{12}x_{45} -
x_{14}x_{25} + x_{24}x_{15}$

\ \ \ \ \ \ \ \ \ \ \ \ \ \ \ \ \ \ \ \ $q_4^+ = uy_4 - x_{12}x_{35} +
x_{13}x_{25} - x_{23}x_{15}$

\ \ \ \ \ \ \ \ \ \ \ \ \ \ \ \ \ \ \ \ $q_5^+ = uy_5 + x_{12}x_{34} -
x_{13}x_{24} + x_{23}x_{14}$

\bigskip

\ \ \ \ \ \ \ \ \ \ \ \ \ \ \ $q_1^- = \ \ \ \ \ \ \ \ \ \ \ \ \ x_{12}y_2 +
 x_{13}y_3 + x_{14}y_4 + x_{15}y_5$

\ \ \ \ \ \ \ \ \ \ \ \ \ \ \ $q_2^- = - x_{12}y_1 \ \ \ \ \ \ \ \ \ \ \ +
x_{23}y_3 + x_{24}y_4 + x_{25}y_5$

\ \ \ \ \ \ \ \ \ \ \ \ \ \ \ $q_3^- = - x_{13}y_1 - x_{23}y_2 \ \ \ \ \ \ \ \ \
\ \ + x_{34}y_4 + x_{35}y_5$

\ \ \ \ \ \ \ \ \ \ \ \ \ \ \ $q_4^- = - x_{14}y_1 - x_{24}y_2 - x_{34}y_3 \ \ \
\ \ \ \ \ \ \ \ + x_{45}y_5$

\ \ \ \ \ \ \ \ \ \ \ \ \ \ \ $q_5^- = - x_{15}y_1 - x_{25}y_2 - x_{35}y_3 -
x_{45}y_4 \ \ \ \ \ \ \ \ \ \ \ $

\end{figure}
\bigskip

An important property of the spinor tenfold
is its self-duality \cite{E}:

\begin{lemma}
The projectively dual variety
$\Sigma^\dual\subset\PP^{15\dual}$, consisting of all the hyperplanes
in $\PP^{15}$ that are tangent to $\Sigma$, is projectively equivalent to 
$\Sigma$.
\end{lemma}

This follows also from
the self-duality of the half-spinor representation of $\Spin$
and the fact that $\Spin$ has only two orbits in $\PP^{15}$:
the spinor tenfold and its complement \cite{I}.

There is an alternative interpretation of the spinor tenfold
$\Sigma$: it is isomorphic to each one of the two families
of 4-dimensional
linear subspaces in a smooth 8-dimensional quadric $Q^8\subset\PP^9$.
In other words, it is a component of the Grassmannian 
$G_q(5,10)=\Sigma^+\sqcup \Sigma^-$ of 
maximal isotropic subspaces of a nondegenerate quadratic form
$q$ in $\CC^{10}$. 
The varieties $\Sigma$,  $\Sigma^\dual$ can be simultaneously identified
with $\Sigma^+$, $\Sigma^-$ respectively in such a way that the duality
between $\Sigma$,  $\Sigma^\dual$ is
given in terms of
certain incidence relations between the four-dimensional
linear subspaces of the quadric $Q^8=\{ q=0 \}$.

Namely,
denote by $\PP^4_c$ the
subspace of $Q^8$ corresponding to a point $c\in \Sigma^\pm$, and
let, for example, $c\in \Sigma^+$. Then we have:

\begin{equation}\label{splus}
\Sigma^+=\{ d\in G_Q(5,10)\mid \dim (\PP^4_c\cap \PP^4_d)
\in \{ 0,2,4\}\}\ ,
\end{equation}

\begin{equation}\label{sminus}
\Sigma^-=\{ d\in G_Q(5,10)\mid \dim (\PP^4_c\cap \PP^4_d)
\in \{ -1,1,3\}\}\ ,
\end{equation}
where the dimension equals $-1$ iff the intersection is empty. Furthermore,
if we denote by $\PP^{14}_w$ the hyperplane in $\PP^{15}$
represented by a point $w\in\check{\PP}^{15}:=(\PP^{15})^\dual$, and by
$H_w$ the corresponding hyperplane section $\PP^{14}_w\cap\Sigma$
of $\Sigma =\Sigma^+$, then for any $w\in \Sigma^-=\Sigma^\dual$, we have

\begin{multline}\label{hw}
H_w = \{ c \in \Sigma^+ : 
                      {\mathbb P}^4_c \cap {\PP}^4_w \not= \empty \} 
        = \{ {\mathbb P}^4 \subset Q : 
 \\  \dim ({\mathbb P}^4 \cap {\mathbb P}^4_w) \mbox{ is odd and } \ge 0 \}
\end{multline}

For future use, we will describe explicitly the identifications
of $\Si$, $\Si^{\dual}$ with $\Si^+$, $\Si^-$. Let $V$ be a
$2\nu$-dimensional $\CC$-vector space ($\nu =5$ in our applications)
with a nondegenerate quadratic
form $q$ and $(\cdot ,\cdot )$ the associated symmetric bilinear form.
Fix a pair of maximal isotropic subspaces $U_0,U_\infty$ of $V$
such that $V=U_0\oplus U_\infty$. The bilinear form $(\cdot ,\cdot )$
identifies $U_0$ with the dual of $U_\infty$. Let $\Si^+$, resp.
$\Si^-$ be the component of $G_q(\nu ,V)$ that contains $U_0$,
resp. $U_\infty$. Let $S=\wedge^{\ccd}U_\infty$ be the exterior
algebra of $U_\infty$ ; it is called the spinor space of
$(V,q)$ and its elements are called spinors. The
even and the odd parts of $S$ 
$$S^+=\wedge^{even}U_\infty,\ 
S^-=\wedge^{odd}U_\infty$$ 
are called half-spinor spaces. To each maximal isotropic
subspace $U\in\Si^+\cup\Si^-$ one can associate a unique,
up to proportionality, nonzero half-spinor
$s_U\in S^+\cup S^-$ such that $\phi_u(s_U)=0$ for all $u\in U$,
where $\phi_u\in\END (S)$ is the Clifford automorphism
of $S$ associated to $u$:
\begin{multline*}
\phi_u(v_1\wedge\ldots\wedge v_k)=\sum_i (-1)^{i-1} (u_0,v_i)v_1
\wedge\ldots\wedge \widehat{v_i}\wedge\ldots\wedge v_k
+u_\infty\wedge v_1\wedge\ldots\wedge v_k\ \\ \mbox{if}\ u=u_0+u_\infty,\
u_0\in U_0,\ u_\infty\in U_\infty .
\end{multline*}

The element $s_U$ is called the pure spinor associated to $U$.
The map $U\mapsto [s_U]\in\PP (S^{\pm})$ is the embedding of
$\Si^{\pm}$ into the projective space $\PP^{2^{\nu -1}-1}$
from which we started our description of the spinor tenfold
(formula (\ref{alt5}), $\nu=5$). The duality between
$S^+$, $S^-$ is given by the so called fundamental
form $\beta$ on $S$, for which $S^+$, $S^-$ are maximal isotropic
$2^{\nu -1}$-dimensional subspaces of $S$:
$$
\beta (\xi ,\xi')=(-1)^{p(p-1)/2}(\xi\wedge\xi')_{top}
$$
where $\deg\xi =p$ and $(s)_{top}$
denotes the $\wedge^\nu U_\infty$-component
of a spinor $s\in \wedge^{\ccd}U_\infty$.

Describe now the spinor embedding in coordinates.
Let $U\in \Si_+$. Then the intersection $U\cap U_\infty$
is always even-dimensional and generically
$U\cap U_\infty =0$. Choose a basis $e_1,\ldots e_\nu$ of 
$U_\infty$ in such a way that $U\cap U_\infty = \langle 
e_{\nu -2k},e_{\nu -2k+1},\ldots ,e_\nu\rangle$.
Let $e_{-1},\ldots , e_{-\nu}$ be the dual basis
of $U_0$. Then $U$ possesses a basis $u_{1},\ldots , u_{\nu}$
of the following form: $u_i=e_{-i}+\sum\limits_{j=1}^{\nu -2k-1} 
a_{ij}e_j$ for $i=1,\ldots ,\nu -2k-1$, and $u_i=e_i$
for $i=\nu -2k,\ldots ,\nu $, where $(a_{ij})$ is a
skew-symmetric matrix of dimension $\nu -2k-1$. The pure
spinor associated to $U$ is given by the following formula:
\begin{equation}\label{sU}
s_U=\exp (\alp )\wedge e_{\nu -2k}\wedge e_{\nu -2k+1}\wedge \ldots 
\wedge e_\nu\ ,\ \ \alp =-\sum\limits_{j=1}^{\nu -2k-1} 
a_{ij}e_i\wedge e_j
\end{equation}
Here the exponential is defined by
$$
\exp\alp =\sum\limits_{i=1}^{[\nu /2]} \frac{\alp\wedge \ldots 
\wedge \alp}{i!}=\sum\limits_{I\subset \{1,\ldots ,\nu -2k-1\}} 
\pf_I(A)e_I
$$
where $e_I=e_{i_1}\wedge \ldots 
\wedge e_{i_p}$ if $I=\{ i_1,\ldots ,i_p\}$
and $\pf_I(A)$ is, up to the sign, the Pfaffian
of the submatrix of $A$ consisting of its rows and columns
with numbers $i_1,\ldots ,i_p$ (so that only even values of $p$
can yield nonzero terms). The coordinates in $\PP^{15}$ used
in formulas (\ref{alt5}), (\ref{equa-sten}) are the ones
corresponding to the basis $1, e_i\wedge e_j,
e_i\wedge e_j\wedge e_k\wedge e_l$ of $S^+=
\wedge^{even}U_\infty$.

\section{Linear sections of the spinor tenfold}
\label{linsect}

Mukai \cite{Mu-1} has observed that a nonsingular section of the
spinor tenfold $\sten\subset\PP^{15}$ by a linear subspace
$\PP^{7+k}$ for $k=-1,0$, resp. 1 is a canonical curve, 
K3 surface, resp. a prime Fano threefold of degree 12. He has proven
that a generic canonical curve of genus 7,
a generic K3 surface of degree 12 and any nonsingular prime Fano threefold 
$X_{12}$ (with Picard group $\ZZ$)
are obtained as linear sections of $\sten$ in a unique
way modulo the action of $\Spin$.

\begin{definition}\label{defdual}
For a given linear section $X$ of $\Sigma$,
we denote its orthogonal linear section by $\check{X}$
and call it the dual of $X$.
In particular, the dual of a Fano linear section $X_{12}$
is a canonical curve $\Gamma =\gseven = \check{X}_{12}$
(the superscript being the genus, and the subscript the degree),
and the dual of a K3 linear section $S$
is another K3 surface $\check{S}$ of degree 12.
\end{definition}

\begin{lemma}\label{nonsing}
Let $\PP^{7+k}$ with $k=-1,0$ or $1$ be a linear subspace in $\PP^{15}$,
transversal to $\Sigma$. Then the orthogonal complement
$(\PP^{7+k})^\perp =\check{\PP}^{7-k}$ is transversal to $\Sigma^\dual$.
Thus there is a natural way to associate to a linear section
of $\Sigma$ which is a Fano threefold,
a K3 surface, resp. a canonical curve, the orthogonal linear section
of $\Sigma^\dual$, which is a canonical curve, a K3 surface,
resp. a Fano threefold of degree $12$.
\end{lemma}

\begin{proof}
Assume that $c\in X=\PP^{7+k}\cap \Sigma$ is a singular point.
We can represent $\PP^{7+k}$ as the intersection of
$8-k$ hyperplanes, so that $X=\Sigma\cap \PP^{14}_{u_0}\cap
\ldots \cap\PP^{14}_{u_{7-k}}$. As $c$ is a singular point,
we can replace the $u_i$ by some linear combinations of them
in such a way that $c\in\Pi^4_{u_0}=\Sing H_{u_0}$.
We can even represent
$X$ as the intersection $\Sigma\cap \PP^{14}_{u_0}\cap
\ldots \cap\PP^{14}_{u_{7-k}}$ with $u_i\in \check{X}$,
since the span of the dual section $\check{X}=\check{\PP}^{7-k}
\cap \Sigma^\dual$ is the whole $\check{\PP}^{7-k}$.

By reflexivity of tangent spaces, $c\in\Pi^4_{u_0}$ implies that
$T_{u_0}\Sigma^\dual\subset H_c=\PP^{14}_c\cap\Sigma^\dual$.
We can complete $c$ to a sequence $c=c_0,\ldots ,
c_{7+k}$ in such a way that $\check{X}=\Sigma^\dual\cap \PP^{14}_{c_0}\cap
\ldots \cap\PP^{14}_{c_{7+k}}$, and the fact that 
$\PP^{14}_{c_0}$ contains the tangent space $T_{u_0}\Sigma^\dual$
implies that $u_0$ is a singular point of $\check{X}$.
We have proven that if $X$ is singular, then $\check{X}$
is. By the symmetry of the roles of $X$ and $\check{X}$, the converse
is also true.
\end{proof}

For future reference, we will cite the following lemma on plane sections
of $\Sigma$. As $\Sigma$ is an intersection of quadrics, every
nonempty section of it by a plane $\PP^2$ is either
a 0-dimensional scheme of length $\leq 4$, or a conic (possibly,
reducible), or a line, or a line plus a point, or the whole
plane. Mukai proves that the case of length 4 is impossible:

\begin{lemma}\label{4-secants}
$\Sigma$ has no $4$-secant $2$-planes.
\end{lemma}

\begin{proof}
This is Proposition 1.16 of \cite{Mu-2}.
\end{proof}

\begin{lemma}\label{singhyp}
For any $w\in \Sigma^-$, the singular locus of $H_w$ is a
projective space $\PP^4$, linearly embedded into $\PP^{15}$,
and it consists only of points of multiplicity $2$.
Denote this $\PP^4$ by $\Pi^4_w$, and its complement
$H_w\setminus\Pi^4_w$ by $U_w$.
Then the linear projection $\pr_w :U_w\lra \PP^9$
with center $\Pi^4_w$ is surjective onto the Grassmannian
$G(2,5)\subset\PP^9$ and induces on $U_w$ the structure
of the universal vector subbundle of $\CC^5\times G(2,5)$ of rank $3$.
\end{lemma}

\begin{proof}
The statement about the multiplicity of $H_w$ at the tangency locus
follows from formulas (\ref{splus})--(\ref{hw}) and
Proposition 2.6 in \cite{Mu-2}, saying that
$\mult_vU_w=\frac{1}{2}(\dim \PP^4_v\cap\PP^4_w+1)$.

The fact that the tangency locus of $\PP^{14}_w$ is a linearly
embedded $\PP^4$ follows from a quite general observation,
which one can refer to as the reflexivity property of
the tangent spaces (see \cite{Kl}): Let $Y\subset
\PP^N$, $Y^*\subset\PP^{N\dual}$ be dual to each other,
$\dim Y=n$, $\dim Y^*=n^*$.
Then for any nonsingular point $[H]\in Y^*$ representing
a hyperplane $H$ in $\PP^N$, the latter is tangent to $Y$ along
the linear subspace $\Pi$ of dimension $N-n^*-1$,
consisting of all the points $[h]\in\PP^N$
such that $T_{[H]}Y^*\subset h$ (a point $[h]\in\PP^N$
represents a hyperplane $h\subset \PP^{N\dual}$). In our case
$N=15$, $n=n^*=10$, so the tangency locus $\Pi$ is $\PP^4$.
Now write down the projection with center $\Pi$ in coordinates.
By homogeneity of $\sten^\dual$, we can choose coordinates
$(u:\hat{X}:\vec{y})$ in such a way that $w=(1:\hat{0}:\vec{0})$
(in dual coordinates), so that the equation of
the hyperplane section is

\medskip 

\centerline{$H_w = \Sigma \cap (u=0)$.} 

\medskip 

In these coordinates, $H_w \subset {\PP}^{14}_w$ 
is defined by the restrictions of the equations (\ref{equa-sten}) 
for $\Sigma \subset {\PP}^{15}$:

\medskip 

\centerline{$\vec{0} = \vec{\pf }(\hat{X}) \ \& \ \hat{X}\vec{y} = 0$.} 

\medskip 

\noindent  
Therefore either $\rk \ \hat{X} = 2$ or $\hat{X} = \hat{0}$, 
and 

\medskip

\centerline{$H_w = U_w \cup {\Pi}^4_w$,}

\medskip 
\noindent 
where   
        $$U_w = \{ (0:\hat{X}:\vec{y}) \in H_c : \rk \hat{X} = 2
	,\  \vec{y}\in\ker\hat{X}\} ,$$
and 
$${\Pi}^4_c = 
   H_c - U_c = \{ (0:\hat{X}:\vec{y}) \in H_c : \rk \hat{X} = 0 \} 
             = \{ (0:\hat{0}:\vec{y}): \vec{y} \in \CC^5 \} .$$
The constraint $\rk \hat{X} = 2$ cuts out exactly the Grassmannian
$G(2,5)$, and $\vec{y}\in\ker\hat{X}$ defines the universal
kernel bundle of rank 3 on it. This proves our
assertion.
\end{proof}

Now let $X=\Sigma\cap\PP^{7+k}$ for $k=-1,0$ or $1$ be a
general linear section of the spinor tenfold, and $\check{X}=
\check{\PP}^{7-k}\cap\Sigma^\dual$ its dual. For any $w\in\check{X}$, let
$\pr_w :U_w\lra G(2,5)$ be the linear projection of Lemma \ref{singhyp}.
We have $X\subset H_w$, and
the nonsingularity of $X$ implies that $X\cap \Pi^4_w=\empty$, that is
$X\subset U_w$.
Let $\pi_w=\pr_w|_X$. Remark that $X$ is a {\em linear} section of $U_w$ and
the fibers of $\pr_w$ are linear subspaces in $U_w$,
so the fibers of $\pi_w$ are also linear subspaces.
They are obviously 0-dimensional if $X$ is
a curve ($k=-1$). As 
$\Pic (X) =\ZZ$ for $k=0,1$, they are 0-dimensional in these cases as well,
and hence $\pi_w$ is a linear isomorphism
of $X$ onto its image in $G(2,5)$. Moreover, $X=\PP^{7+k}\cap U_w$,
hence $\langle X\rangle =\PP^{7+k}$ does not meet $\Pi^4_w$
and $\pi_w(\langle X\rangle )= \langle \pi_w(X)\rangle $
is of dimension $7+k$.

We will now investigate an arbitrary linear embedding
of $X$ into $G(2,5)$.
To this end, we will need Mukai's description of the embedding
of $X$ into the spinor tenfold. Let us forget
that our $X$ is a linear section of $\Sigma$
and construct a spinor embedding of it
in a functorial way. Consider $X$ as a projectively
normal subvariety of some projective space
$\PP^{7+k}$ and denote by $\III_X$ the ideal
sheaf of $X$ in this projective space.
According to
\cite{Mu-1}, \cite{Mu-2}, 
the vector space $V=H^0(\PP^{7+k},\III_X(2))$ is 10-dimensional,
the subspace $U_p=H^0(\PP^{7+k},\III_X(2-2p))\subset V$
is 5-dimensional for any $p\in X$, and this
yields a map $\eta_X:X\lra G(5,V)$, $p\mapsto U_p$.
There is only one quadratic relation between the elements of $V$
(Theorem 4.2, \cite{Mu-2})
providing a quadratic form $q_V$ on $V$, and all the spaces
$U_p$ are maximal isotropic with respect to $q_V$. Thus the
image of $\eta_X$ lies on one of the spinor varieties $\Sigma^{\pm}$
in $G(5,V)$ associated to the quadratic form $q_V$.
Mukai proves (Theorem 0.4, ibid.) that $\eta_X$ is an isomorphism
onto its image. Let us declare
this spinor variety to be $\Sigma^+$, and denote the image of $X$
by $X^+$. Then $\check{X}$ is naturally embedded into
$\Sigma^-$, with image $X^-$, and we can use the incidence
formulas (\ref{splus}), (\ref{sminus}), (\ref{hw}).

\begin{lemma}\label{F8}
Let $X$ be as above, and $i:X\into G(2,5)$ a projective linear embedding,
$U=H^0(\langle G\rangle ,\III_G(2))$ the 
$5$-dimensional space of quadrics passing
through $G$. Assume that the natural map $i^*:U\lra V$ is injective and
that $i^*(U)$ is maximal
isotropic with respect to $q_V$.
Then there exists $w\in X^-$ such that the
map $i\circ\eta_X^{-1}:X^+\into G(2,5)$ and the restriction 
$\pi_w:X^+\lra G(2,5)$ of the
projection $\pr_w$ defined in Lemma \ref{singhyp}
are equivalent under the action of $PGL(5)$ on $G(2,5)$.
\end{lemma}

\begin{proof}
Consider $G=G(2,5)$ in its Pl\"ucker embedding in $\PP^9$ and
identify $U$ with its image in $V$. We have $\dim U=5$ and 
$Z_p=H^0(\langle G\rangle ,\III_G(2-2p))\subset U$ is 
2-dimensional for every $p\in Y$.
This defines a linear isomorphism $\zeta :G\lra G(2,U)$. 
Thus, the original
embedding $i:X\into G(2,5)$ is equivalent to the map
$\zeta_X :=\zeta\circ i:X\into G(2,U)$, sending a point $p\in X$ to
the 2-plane $Z_p=U\cap U_p$. By (\ref{sminus}), 
$w=[U]\in \Sigma^-$, by (\ref{hw}), $X^+\subset H_w$,
and we obtain the linear
projection $\pi_w:X^+\lra G(2,5)$. Let us complete $U=U_\infty$
to a decomposition $V=U_0\oplus U_\infty$ of $V$
into the direct sum of maximal isotropic subspaces. 
Then, as in the proof of Proposition \ref{rho},
(i), $w=s_{U_\infty}\in\Si^-$ and $\pi_w(p)$ is the 
Pl\"ucker image $\xi = \xi_{U_p\cap U_\infty}$
of the 2-plane $Z_p=U\cap U_\infty$.
This ends the proof.

\end{proof}

\begin{lemma} \label{F8-1}
Let $X$ be a nonsingular Fano $3$-fold
($k=1$) or a K3 surface 
of genus $7$ with Picard number $1$ ($k=0$).
Then
for any linear embedding $i:X\into G(2,5)$ such that the map 
$i^*:U=H^0(\langle G\rangle ,\III_G(2))\lra V=
H^0(\langle X\rangle ,\III_X(2))$
is injective,
$U$ is a maximal isotropic subspace
of $V$ with respect to the quadratic
form $q_V$.

\end{lemma}

\begin{proof}
Assume that $U$ is not isotropic. Then $q_V$ defines a 3-dimensional
quadric $Q$ in $\PP (U)$. In the notations from the proof of
Lemma \ref{F8}, the isotropy of the 5-spaces $U_p$
implies that the projective lines $\PP (U\cap U_p)$ ($p\in X$)
are all contained in $Q$. Thus the map $p\mapsto \PP (U\cap U_p)$,
projectively equivalent to $i$, transforms $X$ isomorphically
onto a subvariety of the family of lines $G(1;Q)$
on the 3-dimensional quadric $Q$. 

Let $k=1$, that is $X$ is a Fano threefold.
If $Q$ is nonsingular, $G(1;Q)\simeq\PP^3$,
and this is absurd, as $X\not\simeq\PP^3$. If $Q$ is of rank 4,
then the family of lines on $Q$ has two components, each one of which
is a $\PP^2$-bundle over $\PP^1$; this is absurd because $X$ does
not contain any plane. The cases of smaller rank lead also to contradictions,
hence $U$ is isotropic.

The argument is similar for the case of a K3 surface: if $\rk Q=5$,
then $X\subset \PP^3$, which is absurd, and if $\rk Q=4$, then
$X$ has a pencil of curves defined by the $\PP^2$-bundle over
$\PP^1$, but the generic K3 surface has no pencils of curves.

\end{proof}

Similar statements hold also in the case $k=-1$, but the proof uses
vector bundle techniques and is postponed until Section \ref{map-rho}.

\section{Elliptic quintics and rank-2 vector bundles on $X_{12}$}
\label{ellquint}

Let $X=X_{12}=\PP^8\cap\Sigma$ be a Fano 3-dimensional linear section 
of the spinor tenfold $\Sigma$.
An {\em elliptic quintic} in $X$ is a nonsingular
irreducible curve $C\subset X$ of genus $1$ and of
degree 5. We will also deal with degenerate ``elliptic"
quintics, which we will call just {\em quasi-elliptic quintics}.
A quasi-elliptic quintic is a locally
complete intersection curve $C$ of degree 5 in $X$,
such that $h^0(\OOO_C)=1$ and the canonical sheaf of $C$ is trivial:
$\omega_C=\OOO_C$. A reduced quasi-elliptic quintic will be
called a {\em good quintic}.

\begin{lemma}\label{projnorm}
Let $C\in X$ be a quasi-elliptic quintic. Then
$\langle C\rangle =\PP^4$,
where the angular brackets denote the linear span.
\end{lemma}

\begin{proof}

Assume that $C\subset\PP^3$. Then a general section
of $C$ by a plane $\PP^2\subset \PP^3$ is a 0-dimensional
scheme of length 5. This contradicts
Lemma \ref{4-secants}. Hence $\dim\langle C\rangle \geq 4$.

To prove the opposite inequality, it suffices to show
that $h^0(C, \OOO_C(1))\leq 5$. This follows from
the Serre duality and Riemann--Roch formula.
\end{proof}

Starting from any quasi-elliptic quintic $C\subset X$, one can
construct a rank-2 vector bundle $\EEE$ with Chern
classes $c_1(\EEE )=1$, $c_2(\EEE )=5$.
It is obtained as the middle term of the following
nontrivial extension of $\OOO_X$-modules:
\begin{equation}\label{serre}
0\lra \OOO_X\lra \EEE \lra \III_C(1) \lra 0\; ,
\end{equation}
where $\III_C=\III_{C,X}$ is the ideal sheaf of $C$ in $X$.
One can easily verify (see \cite{MT} for a similar argument)
that, up to isomorphism, there is a unique
nontrivial extension (\ref{serre}), thus $C$ determines
the isomorphism class of $\EEE$. This way of constructing
vector bundles is called Serre's construction. The vector
bundle $\EEE$ has a section $s$ whose scheme of zeros
is exactly $C$. Conversely, for any section $s\in H^0(X,\EEE )$
such that its scheme of zeros $(s)_0$ is of codimension 2,
the vector bundle obtained by Serre's construction from
$(s)_0$ is isomorphic to $\EEE$. As $\det\EEE\simeq\OOO_X(1)$, we
have $\EEE\simeq\EEE^\dual (1)$.

The proofs of the following three lemmas are standard;
see, for example, \cite{MT}, where similar facts
are proved for elliptic quintics in a cubic threefold.

\begin{lemma}\label{elemprop}
For any quasi-elliptic quintic $C\subset X$, the associated vector
bundle $\EEE$ possesses the following properties:

(i) $h^0(\EEE )=5$, $h^i(\EEE (-1))=0\ \forall\ i\in\ZZ$, and
$h^i(\EEE (k))=0\ \forall\ i>0$, $k\geq 0$.

(ii) $\EEE$ is stable and the local dimension of the
moduli space of stable vector bundles at $[\EEE ]$ is at
least 1.

(iii) The scheme of zeros $(s)_0$ of any nonzero section
$s\in H^0(X,\EEE )$ is a quasi-elliptic quintic. 

(iv) If $s,s'$ are two nonproportional sections of $\EEE$,
then $(s)_0\neq (s')_0$. This means that $(s)_0$ and $(s')_0$
are different subschemes of $X$.
\end{lemma}

\begin{lemma}\label{h1ext1}
Let $\EEE $ be a vector bundle as in Lemma \ref{elemprop}, $C$
the scheme of zeros of any nonzero section of $\EEE$,
and $\NNN_{C/X}$ its normal bundle.

Then the following properties are equivalent:

(i) $h^1(\NNN_{C/X})=0$; 

(ii) $h^0(\NNN_{C/X})=5$;

(iii) $h^1(\EEE^\dual \otimes\EEE)=1$;

(iv) $h^2(\EEE^\dual \otimes\EEE)=0$.

If one of the properties (i)--(iv) is verified, then
we have:

(a)
The Hilbert
scheme $\Hilb^{5n}_X$ of subschemes in $X$ with Hilbert polynomial
$\chi (n)=5n$ is smooth and of dimension 5 at
the point $[C]$ representing $C$.

(b) The moduli 
space $M_X=M_X(2;1,5)$ of stable
vector bundles of rank $2$ with Chern classes
$c_1=1,c_2=5$ is smooth and of dimension 1 at the point 
$[\EEE ]$ representing the vector
bundle $\EEE$.

(c) $[\EEE ]$ has a Zariski neighbourhood $U$ in $M_X$
with a universal vector bundle $\EEEE$ over $U\times X$,
and the projective bundle $\PP (\pr_{1*}\EEEE )$ 
is isomorphic to an open subset $H_U$ of $\Hilb^{5n}_X$ via a map 
that can be defined pointwise on the fiber over each point $t\in U$
by $s\mapsto (s)_0$, where 
$s\in (\pr_{1*}\EEEE )_t= H^0(X, \EEE_t)$.

The map $[D]\mapsto [\EEE_D]$ given by Serre's construction
is a smooth morphism from $H_U$ onto $U$
with fiber $\PP^4=\PP H^0(\EEE_D)$.
\end{lemma}

\begin{lemma}\label{ncx}
Let $\EEE $ be a vector bundle as in Lemma \ref{elemprop}, $C$
the scheme of zeros of any nonzero section of $\EEE$,
and $\NNN_{C/X}$ its normal bundle. Then
$h^2(\EEE^\dual \otimes\EEE)\neq 0$ if and only if
$\NNN_{C/X}$ is a decomposable vector bundle on $C$.
In this case $\NNN_{C/X}\simeq \OOO_C\oplus\OOO_C(1)$
and $h^2(\EEE^\dual \otimes\EEE)=h^1(\NNN_{C/X})=1$.
\end{lemma}

Now we will exploit the restrictions of the vector bundles $\EEE$
to the hyperplane sections $H\subset X$, which are K3 surfaces
of genus 7.

\begin{lemma}\label{resH}
Let $\EEE\in M_X(2;1,5)$, $H$ a general hyperplane section of $X$,
and $\EEE_H=\EEE |_H$ the restriction of $\EEE$ to $H$. Then $\EEE_H$
is stable and has the following cohomology:
$h^0(\EEE_H)=5$, 
$h^1(\EEE_H(n))=0$ for all $n\in\ZZ$,
$h^2(\EEE_H(n))=0$ for all $n\geq 0$.
\end{lemma}

\begin{proof}
By \cite{M}, Theorem 3.1, $\EEE_H$ is Mumford--Takemoto semistable
for general $H$. Hence $h^0(\EEE_H(-1))=0$,
and the determinant of $\EEE_H(-1)$ being odd, the semistability
is equivalent to the stability. We have $\chi (\EEE_H)=5$ and
$h^2(\EEE_H)=h^0(\EEE_H(-1))=0$, so $h^0(\EEE_H)\geq 5$. Let $s$
be a section of $\EEE_H$ and $Z=(s)_0$ its scheme of zeros.
It is 0-dimensional, because if it contained a curve from $|\OOO (k)|$,
then $\EEE_H(-k)$ would have nonzero global sections. Thus $s$
defines a Serre triple
\begin{equation}\label{SerreZ}
0\lra \OOO \lra \EEE_H \lra \III_Z(1) \lra 0
\end{equation}
which provides the equivalence
$$
h^0(\EEE_H)=5+k\equi \dim\langle Z\rangle = 3-k,
$$
where the angular brackets denote the linear span.
If we assume that $h^0(\EEE_H)>5$, then $Z\subset\PP^2$,
which contradicts Lemma \ref{4-secants}. Hence $h^0(\EEE_H)=5$,
$h^i(\EEE_H)=0$ for $i>0$. Twisting (\ref{SerreZ}) by $\OOO (n)$,
we deduce the remaining assertions.
\end{proof}

In fact, by the same arguments as above, one proves:

\begin{lemma}\label{onK3}
Let $S$ be any nonsingular surface linear section of $\sten$ by a
subspace $\PP^7$ with Picard
group $\ZZ$. Then a rank-$2$ vector bundle $E$ on $S$ with Chern
classes $c_1=1$, $c_2=5$ is stable if and only if it is obtained
by Serre's construction from a zero-dimensional subscheme
$Z$ of $S$ whose linear span is $\PP^3$. The twists $E(n)$
of any such vector bundle
on $S$ have the same cohomology as $\EEE_H$ in Lemma \ref{resH}.
\end{lemma}

\begin{lemma}\label{fromHtoX}
Let $\EEE\in M_X(2;1,5)$. Then $h^i(\EEE(-1))=0$ for all $i\in\ZZ$, and
$\EEE$ can be obtained
by Serre's construction from a quasi-elliptic
quintic $C\subset X$. Hence $\EEE$ satisfies also the properties
(i), (iii), (iv) of Lemma \ref{elemprop}.
\end{lemma}

\begin{proof}
The exact triple
\begin{equation}\label{resseq}
0\lra \EEE (n-1)\lra \EEE (n)\lra \EEE_H(n)\lra 0
\end{equation}
for generic $H$ together with Lemma \ref{resH} implies
that $h^2(\EEE (n-1))=h^2(\EEE (n))$ and $h^1(\EEE (n-1))\geq h^1(\EEE (n))$
for all $n\geq 0$.
By Kodaira vanishing theorem, $h^2(\EEE (n))=0$ for $n>>0$,
hence $h^2(\EEE (n))=0$ for all $n\geq -1$. Now look at the same
triple for $n=0$. By stability and Serre duality, $h^0(\EEE(-1))=
h^3(\EEE(-1))=0$; and we have just proved that $h^2(\EEE(-1))=0$,
which implies, by Serre duality, that $h^1(\EEE(-1))=0$.
Hence $h^1(\EEE (n))=0$ for all $n\geq 0$. As $\chi (\EEE )=5$,
$h^0(\EEE )=5$. Take any section $s\neq 0$ of $\EEE $. Its
scheme of zeros $C=(s)_0$ is of codimension 2, because
$h^0(\EEE(-1))=0$, so $\EEE$ is obtained by Serre's construction
from~$C$.
\end{proof}

\section{The map $\rho_X:\check{X}\lra M_X$}
\label{map-rho}

Let $X=\PP^{7+k}\cap\Sigma$ for $k=-1,0$ or 1 be a nonsingular linear section 
of the spinor tenfold $\Sigma$, and $\check{X}=
\check{\PP}^{7-k}\cap\Sigma^\dual$ its dual. In the case $k=0$,
assume that $X$ is sufficiently general, so that $\Pic X\simeq \ZZ$.
In the case $k=-1$,
assume that $X$ is a generic curve of genus $7$.

Let $M_X$ be the moduli space $M_X(r;c_1,c_2)$ of stable
vector bundles of rank $r=2$ on
$X$ with Chern classes $c_1=1$, $c_2=5$
in the cases when $X$ is K3 or Fano ($k=0,1$), and the non-abelian
Brill--Noether locus
$W^\alp_{r,K}$ of stable vector bundles on $X$
of rank $r=2$ with canonical determinant $K$ and
with at least $\alp =5$ global sections
in the case when $X$ is a canonical curve ($k=-1$).
In Proposition \ref{rho}, we will construct
a natural morphism $\rho =\rho_X:\check{X}\lra
M_X$.

\begin{lemma}\label{curve-case}
Let $i:\Gamma\lra G(2,5)$ be an embedding of a generic canonical
curve of genus $7$, linear with respect to
the Pl\"ucker coordinates and such that $\langle i(\Gamma )\rangle 
=\PP^6$. Let $\QQQ_{G}$ be the universal quotient rank-$2$
bundle on $G=G(2,5)$ and $\EEE =i^*\QQQ_{G}$.
Then $\EEE$ is stable and $h^0(\EEE )=5$.
\end{lemma}

\begin{proof}
Assume that $h^0(\EEE )<5$. Then there is a section of
$\QQQ_{G}$ vanishing identically on the image of $\Gamma$.
The zero loci of the sections of $\QQQ_{G}$ are the
sub-Grassmannians $G(2,4)$ in $G$, so $i$ embeds $\Gamma$
into some $G(2,4)\subset G$. This is absurd, because the linear
span of $G(2,4)$ is $\PP^5$, but by hypothesis, that of
$i(\Gamma )$ is $\PP^6$. Thus, $h^0(\EEE )\geq 5$
and the restriction map $i^*:H^0(G,\QQQ_{G} )\lra H^0(\Gamma ,\EEE )$
is injective. Denote by $W$ the image of $i^*$. The initial embedding
$i$ is projectively equivalent to the map
$$\phi_W:x\in\Gamma\mapsto 
W_x^\perp =\{ u\in W^\dual \mid 
u(s)=0\ \forall\ s\in W_x\}\in G(2,W^\dual ), $$
where $W_x= \{ s\in W\mid s(x)=0\}$ is of codimension 2
for any $x\in \Gamma$.

Assume that $\EEE$ is non-stable. Then there is
an exact triple
$$
0\lra L_1\lra \EEE \lra L_2\lra 0,
$$
in which $L_1,L_2$ are line bundles, $L_2=\omega_{\Gamma}\otimes L_1^{-1}$
and $\deg L_1\geq 6$. Remark first that the case $h^0 (L_2)=0$
is impossible. Indeed, in this case all the sections of
$\EEE$ are those of the line subbundle $L_1$ and the subspaces
$W_x$ are of codimension 1.

Assume now that $h^0 (L_2)=1$. Then either $W\subset H^0(\Gamma ,L_1 )$
and this brings us to a contradiction as above, or the map
$W\lra H^0(\Gamma ,L_2 )$ is surjective. In the latter case the Pl\"ucker
image of $\phi_W(\Gamma )$ is contained in a linear subspace
$\PP^3$ of $\PP^9$ of the form $\langle e_1\wedge e_5,e_2\wedge e_5,
e_3\wedge e_5,e_4\wedge e_5\rangle$, where $e_1,\ldots ,e_5$
is a basis of $W^\dual$ such that $e_5$ generates the image
of the natural inclusion $H^0(\Gamma ,L_2 )^\dual\lra W^\dual$.
This is absurd, because $\phi_W$ is projectively equivalent to
$i$ and the linear span of $i(\Gamma )$ is $\PP^6$.

Hence $h^0 (L_2)\geq 2$. As $\Gamma$ has no $g^1_4$, we have
$\deg L_2\geq 5$. Hence $\deg L_2=5$ or $6$.

{\em 1$^{st}$ case:} $\deg L_1=\deg L_2=6$. By Riemann--Roch,
$h^0 (L_1)=h^0 (L_2)$. As $h^0(\EEE )$ must be at least
5, we have $h^0 (L_i)\geq 3$, $i=1,2$. Therefore $\Gamma$
has two (possibly coincident) $g_6^2$'s.
By \cite{ACGH}, Theorem V.1.5, the expected dimension of the family $G^r_d$
of linear series $g^r_d$ on a curve of genus $g$
is $\rho_{g,d,r}=g-(r+1)(g-d+r)$, and $G^r_d=\empty$ for
a generic curve of genus $g$
if $\rho_{g,d,r}<0$. Hence a generic $\Gamma$ of genus 7
has no $g_6^2$'s.

{\em 2$^{nd}$ case:} $\deg L_1=7,\deg L_2=5$. If $h^0(L_2)\geq 3$,
then $\Gamma$ has a $g_5^2$, which is impossible by the
same argument as above.
So, $h^0 (L_1)=3, h^0 (L_2)=2$. The Bockstein morphism
$\delta :H^0(L_2)\lra H^1(L_1)$ being given by the cup-product
with the extension class $e\in H^1(L_1\otimes L_2^{-1}) 
=H^0(L_2^{\otimes 2})^\dual$, the vanishing of $\delta$
implies that of $e$. Hence $\EEE =L_1\oplus L_2$.
Let $s_1,s_2,s_3$ be a basis of $H^0(L_1)$
and $t_1,t_2$ that of $H^0(L_2)$. Then $\phi_W$ can be
given in appropriate Pl\"ucker coordinates by
$x\mapsto (s_1(x)e_1+s_2(x)e_2+s_3(x)e_3)\wedge (t_1(x)e_4+t_2(x)e_5)$.
Thus if we fix $t=(t_1:t_2)$ we will get a plane
$\PP^2_t=\langle e_1,e_2,e_3\rangle\wedge (t_1e_4+t_2e_5)$
in $G(2,5)$ which meets $\phi_W(\Gamma )$ in 5 points
in which $(t_1(x):t_2(x))=(t_1:t_2)$. This contradicts
Lemma \ref{4-secants}.

According to \cite{BF}, the non-abelian Brill--Noether loci
$W^\alp_{2,K}$ on a generic curve of genus $g\leq 8$ are empty 
if and only if their expected dimension $d=3g-3-\alp (\alp +1)/2$
is negative. Hence in our case $W^6_{2,K}=\empty$ and
$h^0(\EEE )=5$.
\end{proof}

\begin{proposition}\label{rho}
Denote by $\pr_w :U_w\lra G(2,5)$ for any $w\in \check{X}$ 
the linear projection of Lemma \ref{singhyp}. We have $X\subset U_w=
H_w\setminus\Pi^4_w$.
Let $\pi_w=\pr_w|_X$. It is an isomorphism
of $X$ onto its image in $G(2,5)$.
Define a rank-2 vector bundle $\EEE =\EEE_w$ on $X$
as the pullback of the universal quotient rank-2
bundle on the Grassmannian: $\EEE_w:=\pi_w^*\QQQ_{G(2,5)}$.
Then $\EEE_w $ is stable and
the map $\rho =\rho_X:w\mapsto [\EEE_w ]$
is a morphism from $\check{X}$ to
$M_X$. Any vector bundle $\EEE$ in the image of $\rho$ possesses 
the following properties.

(i) 
$h^0(\EEE )= 5$  and if $k=0$ (resp. $k=1$), then $\EEE$ is
obtained by Serre's construction from a l.~c.~i. 0-dimensional subscheme $Z$
of length 5 such that $\langle Z\rangle=\PP^3$ (resp. from a quasi-elliptic
quintic $C\subset X$).

(ii) $\EEE$ is globally
generated and, for a generic
section $s$ of $\EEE$, $(s)_0$ is a smooth elliptic quintic if $k=1$
and a subset of 5 distinct points if $k=0$.

(iii) Let $k=1$, that is $X$ is a Fano threefold. Then the 
family of singular curves
$(s)_0$ ($s\in H^0(\EEE )$) is a divisor in $\PP H^0(\EEE )$. 
For generic $p\in X$, there are at most three curves $(s)_0$ 
which are singular at $p$.
\end{proposition}

\begin{proof}
In the case $k=-1$, the wanted assertion is an immediate consequence
of Lemma \ref{curve-case}. Consider now the case $k=1$. It implies
easily the one of $k=0$
by taking hyperplane sections.
We will first prove part (i), and 
the stability of $\EEE_w$ will follow from Lemma \ref{elemprop}.

(i) The sections of $\QQQ_{G(2,5)}$ are in a natural one-to-one
correspondence with linear forms $\ell$ on the 5-dimensional
vector space $W$, if we think of $G(2,5)$ as the variety
of 2-planes $\Pi$ in $W$, the fiber $\QQQ_t$ of $\QQQ_{G(2,5)}$
at $t\in G(2,5)$ being just $\Pi^\dual$. Let $s_\ell$
be the section of $\QQQ_{G(2,5)}$ associated to $\ell$;
denote by the same symbol the induced section of $\pr_w^*\QQQ_{G(2,5)}$
and by $s$ its restriction to $X$. 

Let us choose coordinates
$(u:x_{ij}:y_k)$ in $\PP^{15}$ in such a way that the
equations of $\{ s_\ell =0\}$ acquire
a very simple form. First of all, as in the proof of
Lemma \ref{singhyp}, we choose the origin
at $w$, so that $w=(1:\hat{0}:\vec{0})$.
In a coordinate free form, we fix an identification of $\Si$ with $\Si^+$,
as in Section 1, corresponding to a decomposition $V=U_0\oplus U_\infty$
of a 10-dimensional vector space $V$ endowed with a nondegenerate
quadratic form into the direct sum of maximal
isotropic subspaces, and choose $w=s_{U_\infty}\in\Si^-$. 
Then by (\ref{sU}), the $\wedge^2U_\infty$
component of the pure spinor $s_U$ associated to any
maximal isotropic $U\subset V, [U]\in U_w=H_w\setminus \Pi^4_w$,
is just the Pl\"ucker image $\xi = \xi_{U\cap U_\infty}$
of the 2-plane $U\cap U_\infty$; in the notations of (\ref{sU}),
$\xi =e_4\wedge e_5$.
Thus the above 5-space $W$ used for the definition of $G(2,5)$
is naturally
identified with $U_\infty$.
So, if we choose coordinates $(x_1,\ldots ,x_5)$
in $U_\infty$ in such a way that $\ell =x_5$, we obtain
the following equations for the zero locus
of $s_\ell$ in the Pl\"ucker coordinates associated to $(x_1,\ldots ,x_5)$: 
$x_{15}=x_{25} = x_{35} = x_{45} = 0$. 
To these, one should add the equation $u=0$
of $H_w$ and the 10 quadratic ones for $\Sigma =\Si^+$.
Five of the latter ones are trivially satisfied under the above
linear constraints, so finally we obtain the following system of
equations for the closure $Z_\ell$ of
$\{ s_\ell =0\}\subset U_w$ in $\PP^{15}$:

\begin{figure}[!h]

\ \ \ \ \ \ \ \ \ \ \ \ \ \ \ \ \ \ \ \ \ \ \ \ \ 
$u=x_{15} = x_{25} = x_{35} = x_{45} = 0$
\nopagebreak
\medskip

\ \ \ \ \ \ \ \ \ \ \ \ \ \ \ \ \ \  
$q_5^+ = x_{12}x_{34} - x_{13}x_{24} + x_{23}x_{14} = 0$
\nopagebreak
\medskip 

\ \ \ \ \ \ \ \ \ \ \ \ \ \ \ $q_1^- = \ \ \ \ \ \ \ \ \ \ \ \ \ x_{12}y_2 +
 x_{13}y_3 + x_{14}y_4 = 0$
\nopagebreak

\ \ \ \ \ \ \ \ \ \ \ \ \ \ \ $q_2^- = - x_{12}y_1 \ \ \ \ \ \ \ \ \ \ \ +
x_{23}y_3 + x_{24}y_4 = 0$
\nopagebreak

\ \ \ \ \ \ \ \ \ \ \ \ \ \ \ $q_3^- = - x_{13}y_1 - x_{23}y_2 \ \ \ \ \ \ \ \ \
\ \ + x_{34}y_4 = 0$
\nopagebreak

\ \ \ \ \ \ \ \ \ \ \ \ \ \ \ $q_4^- = - x_{14}y_1 - x_{24}y_2 - x_{34}y_3 \ \ \
\ \ \ \ \ \ \ \ = 0$
\end{figure}
\bigskip

The five quadratic equations are just (up to sign)
the quadratic Pfaffians of the skew-symmetric matrix
$$ M=\left[\begin{array}{ccccc} 0&-y_1&y_2&-y_3&y_4\\
y_1&0&x_{34}&x_{24}&x_{23}\\
-y_2&-x_{34}&0&x_{14}&x_{13}\\
y_3&-x_{24}&-x_{14}&0&x_{12}\\
-y_4&-x_{23}&-x_{13}&-x_{12}&0
\end{array}\right] .
$$

By an obvious linear change of variables,
we see that the quadratic Pfaffians of $M$
define the 6-dimensional Grassmannian $G(2,5)$ in the
projective space $\PP^9$ with coordinates
$y_1,y_2,y_3,y_4,x_{12},x_{13},x_{14},x_{23},x_{24},x_{34}$,
and in taking into account the coordinate $y_5$
missing in all the equations, we conclude that $Z_\ell$
is the 7-dimensional cone over $G(2,5)$ with vertex
$(0:\ldots :0:1)\in \PP^{15}$.

It is well known, that the degree of $G(2,5)$ is 5
and that its curve linear sections are 
quintics with trivial canonical bundle \cite{Hu},
so the same property is true for  $Z_\ell$,
if one considers only complete intersection linear sections
(that is, defined by 6 equations).
Adding to the above equations of $Z_\ell$
the 6 linear equations of $X$ in $H_w$, we obtain the
wanted zero locus $(s)_0$ of $s=s_\ell|_X$ as a linear
section of $Z_\ell$ of expected dimension 1. If it is indeed a curve,
we are done. It cannot contain a surface, because the degree of the surface
would not exceed 5, but the Picard
group of $X$ is generated by the class of
hyperplane section which is of degree 12.
Finally, $s$ cannot be identically zero. Indeed, assume
the contrary. The map $\pi_w$ projects $X$ linearly and isomorphically
onto its image $\bar{X}$ in $G(2,5)$, and the fact
that $s\equiv 0$ means that $X$ is contained in the
zero locus of $s_\ell$. The latter is the
Schubert subvariety $\sigma_{11}(L)$, where
$L$ denotes the hyperplane $ \ell =0$ in $V$, that is the
Grassmannian $G(2,4)$ of vector planes contained
in $L$. This is impossible,
because $X$ cannot be represented as a hypersurface in a
4-dimensional quadric. This proves the wanted
assertion about the loci $(s)_0$ and that 
$H^0(G(2,5),\QQQ_{G(2,5)})$ is mapped injectively
into $H^0(X,\EEE_w)$.

(ii) $\QQQ_{G(2,5)}$ is globally generated, hence so is
$\EEE_w=\pi_w^*\QQQ_{G(2,5)}$. The smoothness of the
zero locus of the generic section follows then by Bertini
Theorem.

(iii) Consider $X$ as a subvariety of $G(2,5)$. Let us verify
that for any $p\in X$, there is a Grassmannian $G(2,4)=
\sigma_{11}(L)$ passing through $p$ and such that
its intersection with $X$ is not transversal at $p$, that is,
$\dim T_pX\cap T_pG(2,4)>1$. To this end, choose a basis $e_1,\ldots ,e_5$
of $\CC^5$
in such a way that $p=[e_1\wedge e_2]$. We may assume that
$\sigma_{11}(L_0)\cap X$ is a smooth elliptic quintic,
where $L_0=\langle e_1,e_2,e_3,e_4\rangle$ and 
that $T_pX$ is not contained in the span $\PP^5$
of $\sigma_{11}(L_0)$. Assume also that there is no line
on $\sigma_{11}(L_0)$ through $p$ whose
tangent direction coincides with that of the elliptic quintic; the case
when there is one is treated similarly. 
Under this assumptions we can represent a basis of $T_pX$ in 
the form $(e_1\wedge e_3+e_2\wedge e_4,e_1\wedge e_5+ 
a_{14}e_1\wedge e_4+a_{23}e_2\wedge e_3+a_{24}e_2\wedge e_4,
e_2\wedge e_5+
b_{14}e_1\wedge e_4+b_{23}e_2\wedge e_3+b_{24}e_2\wedge e_4)$,
where $a_{ij},b_{ij}$ are constants.
Any $L$ such that
$p\in \sigma_{11}(L)$ is given by the equation
$\alp_3x_3+\alp_4x_4+\alp_5x_5=0$. Then $T_p\sigma_{11}(L)$
is spanned by four bivectors $e_i\wedge v_j$, $1\leq i,j\leq 2$,
where $(v_1,v_2)$ is a basis of the vector plane 
$\{\alp_3x_3+\alp_4x_4+\alp_5x_5=0\}\subset\langle e_3,e_4,e_5\rangle$.
For example,
if $\alp_5\neq 0$, then one can choose $v_1=-\alp_5e_3+\alp_3e_5$,
$v_2=-\alp_5e_4+\alp_4e_5$.
It is an easy exercise to check that
the $7\times 8$ matrix of components of the vectors
generating $T_pX+T_p\sigma_{11}(L)$ is of rank $<6$ for 
at least one value of $(\alp_3:\alp_4:\alp_5)\in\PP^2$, 
and if the number of such values is finite, then it is
at most three. Since a generic curve $(s)_0\in\PP H^0(\EEE )$ is
smooth, the family of singular ones is at most three-dimensional.
We have seen that the subfamily $Z_p$ of curves $(s)_0\in\PP H^0(\EEE )$
singular at $p$ is nonempty for any $p\in X$, hence, by a dimension
count, $Z_p$ is finite for generic $p$. This ends the proof.
\nopagebreak
\end{proof}

Part (ii) of the Proposition implies the following corollary.

\begin{corollary}
In the case $k=1$, the family of elliptic quintics in $X$ is nonempty.
\end{corollary}

For instance, we have not even verified that the
morphism $\rho :\check{X}\lra
M_X$  is non-constant. This follows from the next lemma.

\begin{proposition} \label{image-open}
The
image of $\rho$ is an irreducible component $M_X^0$ of $M_X$.
\end{proposition}

\begin{proof}
It suffices to prove that the image of $\rho$ is open.
Let $\EEE_0$ be a vector bundle on $X$ in the image
of $\rho$. Then $\EEE_0$ is generated by global sections and
the natural quotient map $H^0(X,\EEE_0)\otimes\OOO_X\lra\EEE_0$
defines a linear embedding of $X$ into the Grassmannian
$G(2,5)$ of 2-dimensional quotients of $\CC^5=H^0(X,\EEE_0)$.
This is an open property and it will be verified for
a vector bundle $\EEE$ in a neighborhood 
of $\EEE_0$ in $M_X$. Also the conditions in the hypotheses of
Lemmas \ref{F8}, \ref{F8-1} are open, so, in the cases when
$X$ is either a K3
surface or a Fano threefold ($k=0$ or $1$),
the embedding
of $X$ into $G(2,5)$ by global sections of $\EEE$ is given,
up to a linear change of coordinates,
by the projection $\pi_w$ for some $w\in \check{X}$, and hence
$\EEE$ is in the image of $\rho$. The case $k=-1$ follows
in the same manner from Lemma \ref{petri-inj} below.
\end{proof}

\begin{lemma}\label{petri-inj} 
Let $X=\Gamma$ be
a generic canonical curve of genus 7.
Under the hypotheses of Lemma \ref{curve-case}, suppose in
addition, as in Lemma \ref{F8-1}, that the map 
$i^*:U=H^0(\langle G\rangle ,\III_G(2))\lra V=H^0(\langle X\rangle ,\III_X(2))$
is injective. Then
$U$ is a maximal isotropic subspace
of $V$ with respect to the quadratic
form $q_V$.
\end{lemma}

\begin{proof}
Assume that $U$ is not isotropic. As in the proof of Lemma
\ref{F8-1}, $q_V$ defines a 3-dimensional
quadric $Q$ in $\PP (U)$ and the isotropy of the 5-spaces $U_p$
implies that the projective lines $\PP (U\cap U_p)$
are all contained in $Q$. Let $\EEE =i^*(\QQQ_{G})$ be
the pullback of the universal quotient rank-2 vector bundle
on $G=G(2,5)$. The fiber of $\EEE$ over $p\in X$
is canonically identified with
the dual of the 2-plane $U\cap U_p$. By Lemma \ref{curve-case},
it is stable and $h^0(X,\EEE )=5$. We can now apply
Proposition 4.1 of \cite{BF}, which yields the injectivity
of the modified Petri map $\Sym^2(H^0(X,\EEE ))\lra H^0(X,\Sym^2(\EEE ))$.
Further, the authors of \cite{BF} prove in the Claim on p. 267
that the injectivity of the modified Petri map is
equivalent to the following property: there is no quadric
$Q$ in $\PP (H^0(X,\EEE )^*)$ containing all the lines $\PP (\EEE_x)$
for $x\in X$. This ends the proof.
\end{proof}

\begin{proposition}\label{irred-curve}
Let $k=-1$, that is $X$ is
a generic canonical curve of genus $7$.
Then $M_X$ is a Fano threefold of genus $7$
and $\rho_X:\check{X}\lra M_X$ is an isomorphism
of Fano threefolds.
\end{proposition}

\begin{proof}
Mukai \cite{Mu-5} has proved that $M_X$ is a Fano threefold of genus 7
with Picard number 1.
By Proposition \ref{image-open}, $\rho_X$ is surjective. It is easy to see that any
non-constant morphism between
two Fano threefolds of genus 7 with Picard number 1
is an isomorphism. Indeed, let $f:X_1\lra X_2$
be such a morphism. The fact that $\Pic X_1\simeq \ZZ$
implies that $f$ is finite of degree $\delta \geq 1$.
Suppose that $\delta > 1$.As $X_1, X_2$ are smooth, the ramification divisors
$\Delta_i\subset X_i$ are smooth surfaces. Let $H_i$
be a hyperplane section of $X_i$. We have, for some positive
integers $d>1,e\geq 1$, the following relations:
$$
f^*H_2\sim dH_1,\ K_{X_i}\sim -H_i,\ f^*\Delta_2\sim\Delta_1,\
K_{X_1}\sim f^*K_{X_2}+(e-1)\Delta_1.
$$
One deduces immediately the relations
$$
\Delta_1\sim \frac{d-1}{e-1}H_1,\ \ \Delta_2\sim 
\frac{d-1}{d}\cdot\frac{e}{e-1}H_2.
$$
As the $\Delta_i$ are integer multiples of $H_i$,
we conclude that $d=e$ and $\Delta_i\sim H_i$ ($i=1,2$).
Hence $\Pic (X_2\setminus \Delta_2)=0$, and this
contradicts the fact that $f$ induces a non-ramified covering of 
$X_2\setminus \Delta_2$. Hence $\delta =1$ and $f$ is an isomorphism.
\end{proof}

\begin{proposition}\label{irredK3}
If $k=0$, that is, $X$ is a K3 surface,
then $M_X=M_X^0$ is a K3 surface; in particular, it is irreducible
and nonsingular. Moreover,
if $X$ is generic, then $\rho_X:\check{X}\lra M_X$
is an isomorphism of K3 surfaces.
\end{proposition}

\begin{proof}
The first assertion is due to Mukai (Proposition 4.4 in \cite{Mu-0}, or
Theorem 6.1.8 in \cite{HL}). The fact that $\rho_X$ is
an isomorphism follows, as above, from Proposition \ref{image-open} and
the fact that there are no surjective morphisms
between K3 surfaces that are not isomorphisms. 
\end{proof}

\begin{corollary}\label{cor-irredK3}
If $k=0$, that is, $X$ is a K3 surface with Picard number 1,
then any $\EEE\in M_X$ is globally generated, and
for generic $s\in H^0(X,\EEE )$, the zero locus $(s)$ is
a set of 5 distinct points which span $\PP^3$.
\end{corollary}
\begin{proof}
By Proposition \ref{irredK3}, $M_X$ is the image of $\rho$ and the
wanted assertions follow from Proposition \ref{rho}.
\end{proof}

\begin{proposition}\label{MXglobgen}
Let $k=1$, that is, $X$ is a generic prime Fano threefold of degree $12$.
Then any $\EEE\in M_X$
is globally generated and the curve $(s)$ is a (smooth)
elliptic quintic for generic $s\in H^0(X,\EEE )$.
\end{proposition}

\begin{proof}
In the setting of Lemma \ref{fromHtoX}, we deduce from the
restriction isomorphism $H^0(X,\EEE )\lra H^0(H,\EEE_H)$ and
from Corollary \ref{cor-irredK3} that for generic $s\in H^0(X,\EEE )$,
the curve $C_s=(s)_0$ is reduced
and may be singular only at a finite set of points $T\subset X$
where $\EEE$ is not generated by global sections. If the restriction
$\EEE_H$ of $\EEE$ to some K3 linear section $H$ of $X$ through
a point $x\in T$ were stable, we would apply the same argument
to see that $\EEE_H$, and hence $\EEE$ itself is globally generated,
which would be a contradiction. Hence $T$
possesses the property that for any $x\in T$ and for
any nonsingular K3 linear section $H$
passing through $x$, $\EEE_H$
is unstable. By Lemma \ref{fromHtoX} and the restriction exact sequence
(\ref{resseq}), $h^0(\EEE_H(-1))=0$, hence $\EEE_H$
is stable for any nonsingular K3 linear section $H$ of $X$ such that
$\Pic (H)\simeq \ZZ$. By Theorem 5.4 of \cite{Mo}, this condition is verified
for a very general $H$ (''very general" means ''in the complement
of a union of countably many Zarisky closed susbsets").
But the proof in \cite{Mo} actually implies
a subtler result in dimension 3: 

\begin{lemma}[Moishezon]\label{moishe}
Let $V$ be a nonsingular projective $3$-fold and $\{ H_t\}_{t\in\PP^1}$
a Lefschetz pencil of hyperplane sections of $V$
such that $h^{2,0}(H_t)>h^{2,0}(V)$ for the nonsingular
members $H_t$ of the pencil. Assume also that the
base locus of the pencil is a nonsingular curve. Then the restriction map
$\Pic (V)\lra \Pic (H_t)$ is an isomorphism for a very
general $t\in\PP^1$.
\end{lemma}

It is easy to see that through any point $x\in X$ passes a Lefschetz
pencil of K3 linear sections satisfying the hypothesis of the lemma.
Hence $\EEE$
is globally generated and $T$ is empty.
\end{proof}

\section{Takeuchi--Iskovskikh--Prokhorov birational maps}
\label{bir-maps}

The general setting for the constructions of the two
birational maps that we will discuss is as follows:
let $X$ be a Fano 3-fold, and $Z\subset X$ a point or a nonsingular curve.
Let $X'$ be the blowup of $Z$. Assume that the linear system
$|-nK_{X'}|$ for some $n>0$ defines a nontrivial birational morphism.
Then it has only finitely many positive-dimensional fibers,
which are curves, there exists a flop 
$X'\dasharrow X^+$
centered on these curves, and the resulting variety $X^+$
possesses an extremal ray that can be contracted down
onto some variety $Y$. Sometimes $Y$ is a 3-fold, and this
is the case studied by Takeuchi--Iskovskikh--Prokhorov.
In what follows, $X=X_{12}$ is a prime Fano threefold of genus 7,
$Z$ will be either a point, or a conic in $X$.

In order to identify the flopping curves we will need several facts
about lines and conics in $X$.

\begin{lemma}\label{lines}
The family of lines on $X$
is parametrized by an equidimensional
reduced curve $\tau (X)$. 
A line $\ell\subset X$ is a regular (resp. singular)
point of $\tau (X)$ if and only
if $\NNN_{C/X}\simeq \OOO_{\PP^1}(-1)\oplus \OOO_{\PP^1}$ (resp.
$\NNN_{C/X}\simeq \OOO_{\PP^1}(-2)\oplus \OOO_{\PP^1}(1)$).
For generic $X$, $\tau (X)$ is a nonsingular curve.

Every line on
$X$ meets only finitely many other lines.
The union of lines $R(X)=\bigcup\limits_{v\in
\tau (X)}\ell_v$ is a surface from the linear system $|\OOO_X(7)|$,
and a generic line meets 8 other lines on $X$. 
\end{lemma}

\begin{proof}
All the assertions, except for the nonsingularity
of $\tau (X)$ for generic $X$, follow from \cite{I2}, Proposition 1,
\cite{I1}, Proposition 2.1, (iv), Proposition 2.4, (iii), Lemma 2.6,
and Theorem 3.1, (vii). Though the statement of the latter Theorem
does not assert that the number of lines meeting the given one is 8,
the proof gives this value (p. 808).

Let us prove the nonsingularity of $\tau (X)$ for generic $X$.
By \cite{RS}, 6.12,
the family of lines $\tau (\Sigma )$ on the spinor tenfold is a
nonsingular irreducible 15-dimensional variety which can be identified with
the family of planes $\PP^2$ contained in the 8-dimensional
quadric $Q^8$. The identification is done as follows:
the line $\ell^\pm_{\PP^2}\subset \Sigma^\pm$ corresponding to a plane
$\PP^2\subset Q^8$ is the pencil $\{ w\in\Sigma^\pm\mid\PP^4_w\supset 
\PP^2\}$. The family of lines contained in the section
of $\Sigma $ by $7$ generic hyperplanes is the scheme of zeros of
a generic section of the rank-14 vector bundle $\VVV =7\QQQ_{\tau (\Sigma )}$
which is the 7-uple direct sum
of the universal rank-2 vector bundle.
As $\QQQ_{\tau (\Sigma )}=\QQQ_{G(2,16)}|_{\tau (\Sigma )}$
and $\QQQ_{G(2,16)}$ is generated by global sections, the
same is true for $\VVV$, so the generic section defines
a nonsingular curve.
\end{proof}

\begin{lemma}\label{conics}
Let $X$ be a generic Fano threefold of degree $12$ with Picard number $1$.
Then the family of conics on $X$ is parametrized by a
generically reduced irreducible 
scheme $\FFF (X)$ of dimension $2$ (the ``Fano surface" of $X$).
A generic conic $C$ is
nonsingular and $\NNN_{C/X}\simeq \OOO_{\PP^1}\oplus \OOO_{\PP^1}$.
The number of conics
passing through a generic point of $X$ is finite.
\end{lemma}

\begin{proof}
By Proposition 4.2.5 and Theorem 4.5.10 of \cite{IP}
the scheme $\FFF (X)$ is generically reduced and
equidimensional of dimension $2$, and the number of conics
passing through a generic point of $X$ is finite. Let us prove
that $\FFF (X)$ is indeed irreducible.

It is obvious that the family of Fano threefold linear sections
of the spinor tenfold $\Sigma$ that contain a fixed conic
in $\Sigma$, if nonempty, is irreducible and of constant
dimension 42. Hence the incidence variety $I$ parametrizing
the pairs $(C,X)$, where $C\subset X$ is a conic and
$X$ is a Fano threefold linear section
of $\Sigma$, is irreducible. The fiber of
its projection $I\lra G(9,16)$ to the second
factor over $X$ is the Fano surface $\FFF (X)$.
By a monodromy argument, to prove the irreducibility
of $\FFF (X)$, it suffices to present a distinguished
component of $\FFF (X)$ for generic $X$. To this end,
we will construct a natural map $\Sym^2\Gamma\lra \FFF (X)$,
where $\Gamma$ is the dual of $X$, and the wanted
distinguished component is the one containing
the image of $\Sym^2\Gamma$.

Identify $\Sigma$ with $\Sigma^+$ and let
$u,v\in\Gamma\subset \Sigma^-$ be two generic points and
$p_{u,v}=\PP^4_u\cap\PP^4_v$. Then
$\dim p_{u,v}\in\{ 0,2,4\}$. The dimension
is definitely not equal to 4. By the description of lines in
$\Sigma$ given in Lemma \ref{lines}, it is not $2$, because
otherwise $\ell^-_{p_{u,v}}$ would be a bisecant line
of $\Gamma$ contained in $\Sigma^-$, but this would contradict
the fact that $\Gamma$ is a linear section of $\Sigma^-$.
So $p_{u,v}$ is a point. The family of all the ${\PP^4}$'s in $Q^8$
in each one of the two components $\Sigma^\pm$ 
passing through a given point $p$
is a 6-dimensional quadric $Q^{6\pm}_p\in\Sigma^\pm$ which can be identified
with an orbit of $\spin\subset\Spin$. 
We can complete
$u=u_1, v=u_2$ to a family of 7 points $u_1,\ldots ,u_7\in\Gamma$
in such a way that $X$ is the intersection of $\Sigma^+$ with 7 hyperplanes
$\PP^{14}_{u_i}$, $i=1,\ldots ,7$. Denote $H_{u_i}$
the hyperplane sections $\Sigma^+\cap \PP^{14}_{u_i}$.

It is easy to see that $Q^{6+}_{p_{u,v}}\subset H_u\cap H_v$.
Indeed, if $\PP^4_+\in Q^{6+}_{p_{u,v}}$, then $u\in \PP^4_u 
\cap \PP^4_+$, hence $\PP^4_u 
\cap \PP^4_+\neq\empty$ and, by (\ref{sminus}), $\PP^4_+\in H_u$.
Similarly, $\PP^4_+\in H_v$. Hence $Q^{6+}_{p_{u,v}}\cap 
H_{u_3}\cap\ldots\cap H_{u_7}\subset X$ is a conic.

We have constructed a rational map $f:\Gamma^{(2)}\lra \FFF (X)$,
and this ends the proof.

\end{proof}

We will use the symbol ${\mathcal C}^g_d[k]_Z$
to denote the family of all the connected curves
of genus $g$, degree $d$, meeting $k$
times a given subvariety $Z$. More precisely, let
$Z\subset X$ be a nonsingular curve (resp. a point).
Then ${\mathcal C}^g_d[k]_Z$ is the closure in
the Chow variety of $X$ of the family of reduced connected curves
$C$ of degree $d$ such that 
length$\: (\OOO_X/(\III_C+\III_Z))=k$ (resp. $\mult_ZC=k$) and
$p_a(\tilde{C})=g$, where $\tilde{C}$ is the proper transform
of $C$ in the blowup of $Z$ in $X$.

\subsection*{Birational isomorphisms $\Phi_p:X_{12} \rightarrow Y_5$}

These are birational isomorphisms from a given
variety $X=X_{12}$
to a Del Pezzo threefold $Y=Y_5$ of degree $5$ parametrized by
a sufficiently general point $p$ in $X$. The threefold $Y_5$
is defined as a smooth section of the
Grassmannian $G(2,5)\subset\PP^9$ by a subspace
$\PP^6$. It is a Fano threefold of index two,
that is $-K_Y$ is twice the class of hyperplane section
of $Y$.

\begin{theorem}
\label{bir1}
Let $X = X_{12} \subset {\PP}^8$ 
be a smooth anticanonically ebbedded Fano $3$-fold of
index $1$ and of degree $12$, and let $p \in X$ be a sufficiently
general point of $X$. ``Sufficiently general" here means that
it does not lie on a line in $X$  
and satisfies the condition (4.1)(**) from \cite{IP}.
Then the following assertions hold: 

{\bf (a)} \ The non-complete linear system $\mid {\cal O}_X(3 - 7p) \mid$
defines a birational isomorphism ${\Phi_p=\phi}:  X \rightarrow Y$, where $Y = Y_5$ is
the Fano 3-fold of index 2 and of degree 5.

{\bf (b)} \ The one-dimensional family ${\cal C}^0_7[3]_p$ of curves $C \subset
X$ of degree $7$, of genus $0$ and such that ${mult}_p(C) \ge 3$ sweeps out the
unique effective divisor $M = M_p$ in the linear system $\mid {\mathcal
O}_X(5 - 12p) \mid$.

{\bf (c)} \ The birational map ${\phi}:  X \rightarrow Y$ can be represented as
a product $\phi = {\tau} \circ {\kappa} \circ {\sigma}^{-1}$, where ${\sigma}:X'
\rightarrow X$ is the blowup of $p \in X$, ${\kappa}:X \rightarrow X^+$ is
a flop over the double projection $X''$ of $X$ from $p$, and ${\tau}:X^+ \rightarrow Y$
is a blow-down of the proper image $M^+ \subset X^+$ of $M$ onto a canonical
curve ${\Gamma} \subset Y \subset {\bf P}^6$ of degree $12$ (and of genus $7$).
The extremal curves $C^+ \subset X^+$ contracted by $\tau$ 
are the strict transforms of the curves $C
\in {\mathcal C}^0_7[3]_p$.  

{\bf (d)} \ The birational map ${\psi} = {\phi}^{-1}:  Y \rightarrow X$ is
defined by the non-complete linear system $\mid {\mathcal O}_Y(12 - 7{\Gamma})
\mid$.

{\bf (e)} \ The two-dimensional family ${\mathcal C}^0_7[12]_{\Gamma}$ of curves $C
\subset Y$ of degree $7$ and genus $0$ meeting $\Gamma$ $12$
times sweeps out the unique effective divisor $N = N_{\Gamma}$ in
the linear system $\mid {\mathcal O}_Y(5 - 3{\Gamma}) \mid$.  The
proper transform $N' \subset X'$ of $N$ coincides with the exceptional divisor
${\sigma}^{-1}(p) \cong {\bf P}^2$ of $\sigma$.

The extremal curves $C' \subset X'$ contracted by $\sigma$ 
are the strict transforms of the curves $C \in
{\mathcal C}^0_7[12]_{\Gamma}$.

{\bf (f)} \  Let $q_1,...,q_e \subset X$ be all the conics on $X$
which pass through $p$, 
$q_1',...,q_e' \subset X'$ the proper transforms of $q_1,...,q_e$ 
on $X'$. Let $l_1,...,l_{e'}$ be all the bisecant lines to $\Gamma$ 
which lie on $Y$, 
and $l_1^+,...,l_{'e}^+$ the proper transforms of $l_1,...,l_{e'}$ 
on $X^+$. Then for generic $p\in X$ we have $e=e'=24$,
$q_1',...,q_e' \subset M'$ (resp. $l_1^+,...,l_e^+\subset N^+$)
are all the flopping curves of $\kappa$ (resp. $\kappa^{-1}$), and
$\kappa$ transforms $q'_i$ into $l^+_i$ 
for an appropriate ordering of the $l_i$ ($i=1,\ldots ,e$). The map
$\pi'_{2p}$ (resp. $\pi^+_{\Gamma}$) is a small birational morphism contracting
the curves $q_i'$ (resp. $l^+_i$) into isolated ordinary double points.
It is given by the linear system $|K_{X'}|=\sigma^*|\OOO_X(1)-2p|$
(resp. $|K_{X^+}|=\tau^*|\OOO_{Y}(2)-\Gamma |$) and its image 
$\overline{X} = \overline{Y}$ is a quartic threefold in $\PP^4$
with only $e=24$ singular points $\pi'_{2p}(q'_i)$
(resp. $\pi^+_{\Gamma}(l^+_i)$), $i=1,\ldots ,e$. 
\end{theorem}

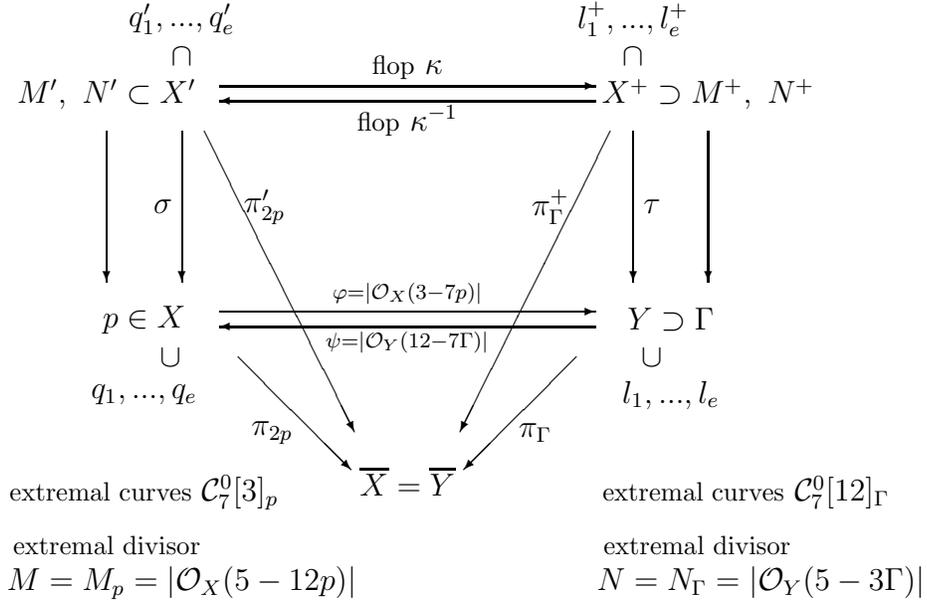
\begin{figure}[!h]
\begin{picture}(120,75)

\put(30,61){\vector(1,0){50}} 
\put(80,59){\vector(-1,0){50}}
\put(15,55){\vector(0,-1){20}}
\put(25,55){\vector(0,-1){20}}
\put(85,55){\vector(0,-1){20}}
\put(95,55){\vector(0,-1){20}}
\put(28,55){\vector(1,-2){20}}
\put(82,55){\vector(-1,-2){20}}
\put(30,31){\vector(1,0){50}}
\put(80,29){\vector(-1,0){50}}
\put(32.5,25){\vector(1,-1){15}}
\put(77.5,25){\vector(-1,-1){15}}
\put(25,70){\makebox(0,0){$q_1',...,q_e'$}}
\put(85,70){\makebox(0,0){$l_1^+,...,l_e^+$}}
\put(55,63.5){\makebox(0,0){{\footnotesize flop} $\kappa$}}
\put(55,56.5){\makebox(0,0){{\footnotesize flop} $\kappa^{-1}$}}
\put(25,65){\makebox(0,0){$\cap$}}
\put(85,65){\makebox(0,0){$\cap$}}
\put(15,60){\makebox(0,0){$M', \ N' \subset X'$}}
\put(95,60){\makebox(0,0){$X^+ \supset M^+, \ N^+$}}
\put(22.5,45){\makebox(0,0){${\sigma}$}}
\put(87.5,45){\makebox(0,0){${\tau}$}}
\put(36,45){\makebox(0,0){$\pi'_{2p}$}}
\put(74,45){\makebox(0,0){${\pi}^+_{\Gamma}$}}
\put(55,33.5){\makebox(0,0){${\scriptstyle \varphi = |{\mathcal O}_X(3-7p)|}$}}
\put(55,27){\makebox(0,0){${\scriptstyle \psi = |{\mathcal O}_Y(12-7\Gamma )|}$}}
\put(20,30){\makebox(0,0){$p \in X$}}
\put(90,30){\makebox(0,0){$Y \supset \Gamma$}}
\put(23.5,25){\makebox(0,0){$\cup$}}
\put(87.5,25){\makebox(0,0){$\cup$}}
\put(20,20){\makebox(0,0){$q_1,...,q_e$}}
\put(90,20){\makebox(0,0){$l_1,...,l_e$}}
\put(37,15){\makebox(0,0){${\pi}_{2p}$}}
\put(72,15){\makebox(0,0){${\pi}_{\Gamma}$}}
\put(55,8.5){\makebox(0,0){$\overline{X} = \overline{Y}$}}

\put(20,7){\makebox(0,0)
      {{\footnotesize extremal curves} 
      ${\mathcal C}^0_7[3]_p$}}
\put(100,7){\makebox(0,0)
      {{\footnotesize extremal curves} 
      ${\mathcal C}^0_7[12]_{\Gamma}$}}
\put(15,0){\makebox(0,0)
      {{\footnotesize extremal divisor}}} 
\put(25,-5){\makebox(0,0)
      {$M = M_p = |{\mathcal O}_X(5 - 12p)|$}} 
\put(93.5,0){\makebox(0,0)
      {{\footnotesize extremal divisor}}}
\put(102,-5){\makebox(0,0)
      {$N = N_{\Gamma} = |{\mathcal O}_Y(5 - 3{\Gamma})|$}} 

\end{picture}

\bigskip
\bigskip
\bigskip
\bigskip
\caption{The birational isomorphism between $X_{12}$ and the Del Pezzo
3-fold $Y_5 = G(2,5) \cap {\bf P}^6$ defined 
by a point $p \in X_{12}$.}

\bigskip\label{fig-1}
\end{figure}

The statement of the theorem is illustrated by Diagram \ref{fig-1}.

\noindent

\begin{proof}
The parts (a)-(e) are essentially an expanded version of Theorem 4.5.8, (ii)
from Iskovskikh--Prokhorov \cite{IP}. According to loc. cit.,
the linear system  of the divisor $N$ is nothing but the proper transform
of $|-5K_{X^+}-2M^+|$, which coincides with our
$\mid {\mathcal O}_X(5 - 12p) \mid$. Other values of $\alp, \beta$ in various
linear systems $|\OOO_X(\alp -\beta p)|$ or
$|\OOO_Y(\alp -\beta \Gamma )|$ are extracted from the preliminary
material of Chapter 4 in \cite{IP} (Lemmas 4.1.1--4.1.9).

Now turn to the proof of (f). By \cite{IP}, Lemma 4.1.1,
the small morphisms  $\pi'_{2p}$, $\pi^+_{\Gamma}$ are given by $n$-th
multiples of the anticanonical linear system for some $n>0$.
By ibid., the proof of
Proposition 4.5.10, $n=1$ and the image is a quartic threefold.
By Lemma 4.5.1, the image $F$ of any flopping curve of $\kappa$ in $X$
belongs to one of the three
types: a conic passing through $p$, a quartic with a double point
at $p$, or a sextic with a triple point at $p$. But the
second and the third types are impossible, because such
curves $F$ span $\PP^3$, and hence $F$ has a 3-dimensional
family of 4-secant planes. It defines a 3-dimensional family of
conics in $X$ by Lemma \ref{4-secants}, which contradicts
the fact that the family of conics on $X$ is 2-dimensional.
Hence the only flopping curves in $X$ are the conics passing through $p$.

By Lemma \ref{conics}, all the conics
passing through $p$ have normal bundle
$\NNN_{q_i/X}\simeq \OOO_{\PP^1}\oplus \OOO_{\PP^1}$. They
also have different tangent directions at $p$ (see Remark \ref{*-**}
below), so the curves $q_i'$ are disjoint. Hence
$\NNN_{q_i'/X'}\simeq \OOO_{\PP^1}(-1)\oplus \OOO_{\PP^1}(-1)$
and $\kappa$ consists in blowing up all the curves $q_i$
and blowing down the obtained exceptional quadrics $\PP^1\times \PP^1$
along the other ruling. For such a one-storey flop, its defect
$(N')^3-(N^+)^3$ is just the number $e$ of flopping curves.
It can be determined from the second formula
(2.6.4) of \cite{Tak}; one obtains $e=24$.

To identify the images of $q_i'$ in $Y$, remark that the images $q_i^+$ of
$q_i'$ under $\kappa$ are contained in $N^+$ and satisfy
$(q_i^+)^2_{N^+}=(q_i\cdot N^+)_X=-1$. As $N^+\sim 
5\tau^*H-3M^+$, where $H$ is the hyperplane section of $Y$,
we obtain the equation $5\deg \tau (q_i^+)-3\length (\tau (q_i^+)\cap\Gamma ) 
=-1$. Further, by the general definition of a flop, all the flopping
curves have zero intersection number with the canonical class,
so $(q_i^+\cdot K_{X^+})=0$, which gives the equation
$2\deg \tau (q_i^+)-\length (\tau (q_i^+)\cap\Gamma )=0$.
Thus $\deg \tau (q_i^+)=1$, $\length (\tau (q_i^+)\cap\Gamma )=2$,
that is, $\tau (q_i^+)$ is a bisecant line to $\Gamma$.
This ends the proof of the Theorem.
\end{proof}

The first application of this theorem is the following lemma.

\begin{lemma}\label{quintp}
Let $p$ be a generic point of $X$. Then all the curves
$C\in {\mathcal C}^0_5[2]_p$, except for finitely many
of them, are
irreducible good quintics with only one double point (node or cusp) at $p$,
$\dim {\mathcal C}^0_5[2]_p=1$, and there is a rational dominant
map of curves $\pi :{\mathcal C}^0_5[2]_p\lra \Gamma_{12}^7$
of degree $\leq 3$, where $\Gamma_{12}^7$ is the genus-$7$
curve from the statement of Theorem~\ref{bir1}.
\end{lemma}

\begin{proof}
The map $\Phi_p$ transforms the curves $C$ from ${\mathcal C}^0_5[2]_p$
into lines meeting $\Gamma$ one time, so all of them are
rational curves with a node or a cusp, except for finitely many ones
which are transforms of lines in $Y_5$ meeting the flopping curves. 
The structure of the family of lines on $Y_5$
is well known (see e. g. \cite{Ili}, \cite{FN}). There is only
one variety $Y_5$ up to projective equivalence, and it possesses a $SL(2)$-action
with three orbits: $C_6^0$, $S_{10}\setminus C_6^0$ and
$Y_5\setminus S_{10}$, where $C_6^0$ is a rational normal sextic
curve and $S_{10}$ its tangential scroll (the subscript
indicates the degree). The family of lines on $Y_5$
is parametrized by $\PP^2$ and through any point of the $k$-th
orbit for $k=1,2,3$ (in the above order) pass exactly
$k$ lines. Define $\pi$ by assigning to each
$C\in {\mathcal C}^0_5[2]_p$ the point of intersection
of $\Phi_p(C)$ with $\Gamma_{12}^7$.
\end{proof}

\subsection*{Birational isomorphisms $\Psi_q$ with a three-dimensional quadric}

This is the family of birational maps $\Psi_q:X_{12}\lra Q_2$
to the three-dimensional quadric $Q_2\subset\PP^4$
parametrized by a sufficiently general conic $q\subset X$.

\begin{theorem}
\label{bir2}
Let $X = X_{12} \subset {\PP}^8$ 
be a smooth anticanonically ebbedded Fano $3$-fold of
index $1$ and of degree $12$, and let $q\subset X$ be a sufficiently
general conic in $X$. ``Sufficiently general" here means that it
satisfies the conditions (4.1) (*), (**) from \cite{IP}.
Then the following assertions hold: 
\smallskip 

{\bf (a)} \ The non-complete linear system $\mid {\cal O}_X(2 - 3q) \mid$
defines a birational isomorphism ${\Psi_q=\phi}:  X \rightarrow Q$, where $Q = Q_2$ is
a smooth three-dimensional quadric.

\smallskip 

{\bf (b)} \ The one-dimensional family ${\cal C}^0_3[2]_q$ of twisted cubics $C
\subset X$ meeting twice the conic $q$, sweeps out the unique effective
divisor $M = M_q$ in the non-complete linear system 
$\mid {\cal O}_X(5 - 8q) \mid$.

\smallskip 

{\bf (c)} \ The birational map ${\phi}:  X \rightarrow Y$ can be represented as
a product $\phi = {\tau} \circ {\kappa} \circ {\sigma}^{-1}$, where ${\sigma}:X'
\rightarrow X$ is the blowup along $q \subset X$, ${\kappa}:X \rightarrow
X^+$ is a flop over the projection $X''$ of $X$ from $q$, and ${\tau}:X^+
\rightarrow Y$ is the blow-down of the proper image $M^+ \subset X^+$ of $M$ onto
a curve ${\Gamma} \subset Q$ of degree $10$ and of genus $7$.

The extremal curves $C^+ \subset X^+$ contracted by $\tau$ 
are the strict transforms of the curves $C
\in {\cal C}^0_3[2]_q$. 

\smallskip 

{\bf (d)} \ The birational map ${\psi} = {\phi}^{-1}:  Q \rightarrow X$ is
defined by the linear system $\mid {\cal O}_Q(8 - 3{\Gamma}) \mid$.

\smallskip 

{\bf (e)}  \ The one-dimensional family ${\cal C}^0_3[8]_{\Gamma}$ of twisted
cubics on $Q$ intersecting the curve $\Gamma$ in a divisor of degree $8$, sweeps
out the unique effective divisor $N = N_{\Gamma}$ in the linear
system $\mid {\cal O}_Q(5 - 2{\Gamma}) \mid$. The proper transform $N' \subset X'$
of $N$ coincides with the exceptional divisor ${\sigma}^{-1}(q) \cong {\bf P}^1
\times {\bf P}^1$ of $\sigma$.

The extremal curves $C' \subset X'$ contracted by $\sigma$ 
are the strict transforms of the curves $C \in
{\cal C}^0_3[8]_{\Gamma}$.
\end{theorem}

The statement of the theorem is illustrated
by Diagram \ref{fig2}.

\begin{figure}[h]

\noindent
\begin{picture}(120,40)(20,0)\thicklines \linethickness{0.05mm}
\put(75,33){\makebox(0,0){${\scriptstyle {\phi}\ =\ |O_X(2-3q)|}$}}
\put(40,28){\makebox(0,0){{ conic} $q \subset X_{12}$}}
\put(55,29){\vector(1,0){40}}
\put(105,28){\makebox(0,0){$Q_2 \supset
{\Gamma}^7_{10}$}}
\put(95,27){\vector(-1,0){40}}
\put(75,23){\makebox(0,0){${\scriptstyle {\psi}\ =\ |O_Q(8-3{\Gamma})|}$}}
\put(40,17){\makebox(0,0)
           {{\footnotesize extremal curves} ${\cal C}^0_3[2]_q$}}
\put(115,17){\makebox(0,0)
            {{\footnotesize extremal curves} ${\cal C}^0_3[8]_{\Gamma}$}}
\put(35,12){\makebox(0,0)
           {{\footnotesize extremal divisor}}}
\put(40,7){\makebox(0,0){$M_q = |O_X(5 - 8q)|$}}
\put(110,12){\makebox(0,0)
            {\footnotesize extremal divisor}}
\put(115,7){\makebox(0,0){$N_{\Gamma} = |O_Y(5 - 2{\Gamma})|$}}
\end{picture}

\bigskip

\caption{The birational isomorphism between $X_{12}$ and 
the quadric
3-fold  $Q_2 \subset {\bf P}^4$
defined by a conic $q \subset X_{12}$.}
\label{fig2}
\end{figure}
\bigskip

\begin{proof}
See Iskovskikh--Prokhorov \cite{IP}, Theorem 4.4.11, (iii).
\end{proof}
\begin{remark}\label{*-**}
According to \cite{IP}, p. 66 (Remark) and Proposition 4.4.1,
the conditions (4.1)(*)-(**) will be satisfied for a smooth 
conic $q \subset X$ if:
 
a) there is a finite number of lines on $X$ which intersect $q$;

b) there is a finite number of bisecant conics to $q$ on $X$. 

It is easy to see that the property b) is automatically verified and in fact
even a stronger condition is satisfied: there are no pairs of bisecant
conics in $X$.  Indeed, the union
of two bisecant conics has a 3-dimensional family of 4-secant
planes $\PP^2$, each of which determines a conic contained in $X$
by Lemma \ref{4-secants}. This is absurd, because the family of
conics in $X$ is 2-dimensional.
It seems very likely that the property a) is also verified
for any smooth conic on $X$. 
\end{remark}

The second birational isomorphism also provides examples of
good quintics in $X$:

\begin{lemma}\label{quintq}
Let $q$ be a general conic in $X$. Then all but a finite number
of curves of the form $q\cup C_3^0$ with $C_3^0\in {\cal C}^0_3[2]_q$
are good quintics. The family of good quintics of this
form is birationally parametrized by an open set of the genus-$7$ curve $\Gamma_{10}^7$
from Theorem \ref{bir2}.
\end{lemma}

\begin{proof}
By part (c) of the above theorem, the strict transforms
of curves from ${\cal C}^0_3[2]_q$ with respect to
$\sigma$ are the extremal curves of the contraction
$\tau$, that is, they form the family of $\PP^1$'s in $M^+$
contracted to the points of $\Gamma$.
\end{proof}

\subsection*{Three curves of genus 7 associated to $X$}

So far, we have associated to the Fano threefold $X=X_{12}$
three curves of genus $7$: the dual linear section $\check{X}$,
the curve $\Gamma_{12}^7$ from Theorem \ref{bir1}, and
$\Gamma_{10}^7$, the one from Theorem \ref{bir2}.
Denote them 
by different symbols, say, $\Gamma$, $\Gamma '$ and $\Gamma ''$
respectively. We will see right now that $\Gamma '$, $\Gamma ''$
are isomorphic; the isomorphism with $\Gamma$ will be deduced
later from the birational isomorphisms of $\Gamma$, $\Gamma ''$ 
to the same component of the moduli space $M_X=M_X(2;1,5)$
of vector bundles on $X$.

\begin{lemma}
The intermediate Jacobian $J(X)$ is isomorphic to the Jacobians
of both curves $\Gamma '$, $\Gamma ''$, so by Torelli Theorem,
$\Gamma '\simeq \Gamma ''$.
\end{lemma}

\begin{proof}
By \cite{CG}, the intermediate Jacobian of the
blowup $\tilde{Y}$ of a smooth curve $D$
in a smooth projective 3-fold $Y$ is isomorphic,
as a p.~p.~a.~v., to the product $J(Y)\times
J(D)$. As $J(Q_2)=J(Y_5)=0$ and the birational isomorphisms $\Phi_p$,
$\Psi_q$ are obtained by blowing up or down
only one nonrational curve and some sets of rational
curves, we have $J(X)\simeq J(\Gamma ')\simeq
J(\Gamma '')$, so the assertion follows by Torelli Theorem.
\end{proof}

The good quintics from Lemmas \ref{quintp}, \ref{quintq}
define some vector bundles of the type considered in
Section \ref{ellquint}. We will see that they fill an
irreducible component of the moduli space.

\begin{lemma}
Let $M_X=M_X(2;1,5)$ be as in Section \ref{ellquint}. Let $p\in X$
be a generic point and $q\subset X$ a generic conic.
Let $M_X(p)$, resp. $M_X(q)$ be the closure in
$M_X$ of the locus of vector bundles $\EEE_C$ on $X$
obtained by Serre's construction from generic good quintics
$C\in {\mathcal C}^0_5[2]_p$, resp. $C= q\cup C_3^0$ 
with $C_3^0\in {\cal C}^0_3[2]_q$. Then
$M_X(p)$, $M_X(q)$ coincide with a $1$-dimensional
irreducible component $M_X^1$ of $M_X$ which does not
depend on $p,q$. Moreover, a generic vector bundle from $M_X^1$
can be obtained by Serre's construction from a smooth elliptic quintic,
and $(M_X^1)_{\red}$ is birational to $\Gamma ''$.
\end{lemma}

\begin{proof}

Denote by $M_X^1$ the closure of the union of all the $M_X(p)$, where $p$
runs over the open subset $U\in X$ for which the conditions of Theorem
\ref{bir1} are satisfied. A generic $\EEE_0$ in any component of $M_X^1$
is obtained by Serre's construction from a good quintic, so
by Lemma \ref{elemprop},  
the family of curves $(s)_0$
of zeros of all the sections $s$ of $\EEE_0$ is $\PP^4=\PP H^0(\EEE_0 )$.
By Proposition \ref{MXglobgen}, the generic curve $(s)_0$ is nonsingular.
As the locus of points of $\PP^4$ corresponding to singular quintics
is non-empty, it is a hypersurface in $\PP^4$, hence it is 3-dimensional.
The family $\bigcup_{p\in U}{\mathcal C}^0_5[2]_p$ being 4-dimensional,
the subfamily of these quintics contained in 
$\PP H^0(\EEE_0 )$ for generic $\EEE_0$ is 3-dimensional; if we fix now
such a generic $\EEE_0$, then
the subfamily of quintics
with only one node at $p$ is nonempty for any $p\in U$ and forms a curve for
generic $p$. In particular, $\EEE_0\in M_X(p)$ for any $p\in U$.
Hence $M_X^1=M_X(p)$ for generic $p$. Further, $M_X^1$ is a union 
of components of $M_X$. Indeed, the same argument as above
shows that any $\EEE\in M_X$
close to a generic $\EEE_0\in M_X^1$ also belongs to $M_X^1$.
Finally, by Lemma \ref{quintp}, $M_X^1$ is birational
to the family of lines in $Y_5$
meeting $\Gamma '$, so it dominates $\Gamma '$ with degree
at most 3.

Let us show now that $M_X(q)$ for a generic conic $q$ is contained in $M_X^1$.
Indeed, the structure of $\Psi_q$ as described in Theorem \ref{bir2}
implies easily that the two points $p,p'\in q$
of the intersection $C_3^0\cap q$
are movable on $q$, when $C_3^0$ runs over ${\cal C}^0_3[2]_q$, so for
generic $q$ and generic $C_3^0$, any of them is a generic point of $X$.
Thus the curve $C=C_3^0\cup q$
for generic $C_3^0\in{\cal C}^0_3[2]_q$ is an element
of ${\mathcal C}^0_5[2]_p$, where $p$ satisfies the hypotheses
of Theorem \ref{bir1}. 
Let $\EEE\in  M_X(q)$ be the vector bundle associated to $C$.
By Lemma \ref{elemprop},  
the family of curves $(s)_0$
of zeros of all the sections $s$ of $\EEE$ is $\PP^4=\PP H^0(\EEE )$.
The family of curves of type $C_3^0\cup q$ ($C_3^0\in{\cal C}^0_3[2]_q$) being
3-dimensional, $C$ deforms inside $\PP H^0(\EEE )$ into a curve
with only one node, hence $\EEE$ is in the closure of the union of
$M_X(x)$, $x\in X$ generic. Hence $M_X(q)\subset M_X^1$ is an irreducible
component, birational to $\Gamma ''$.

It remains to see that $M_X^1$ is irreducible. Let 
$\EEE\in M_X(p)$ for generic $p$. Then $\PP H^0(\EEE )$
contains a 3-dimensional family of curves $D$ with two nodes $x,x'$ for
which each one of the points $x,x'$ runs over a dense 
subset of $X$. As $p_a(D)=1$, the generic $D$ is reducible.
There is no line through a generic point of $X$, so the only possible
type of decomposition of $D$ is a conic plus a cubic.
The cubic is rational and nonsingular, because $X$ is an intersection of quadrics.
Hence $\EEE\in M_X(q)$ and $M_X^1=M_X(q)$ for some $q$.
\end{proof}

As $\Gamma '$, $\Gamma ''$
are curves of the same genus, $\Gamma ''$ does not have
a dominant rational map to $\Gamma '$ of degree
$>1$, so we can deduce from the proof of
the lemma the following corollary:

\begin{corollary}\label{threequintics}
Let $p\in X$
be a generic point.
Let $\EEE_C$ be the vector bundle on $X$
obtained by Serre's construction from a generic good quintic
$C\in {\mathcal C}^0_5[2]_p$ and $\ell =\Phi_p(C)$ the line
in $Y_5$ meeting $\Gamma '$ at a point $u$.
Then there are at most three distinct lines $\ell_i$
in $Y_5$ meeting $\Gamma '$ at $u$ ($\ell_1=\ell$),
and the vector bundles obtained by Serre's construction
from the curves $C_i:=\Phi_p^{-1}(\ell_i)$ are isomorphic to
$\EEE$. Equivalently, the family $\PP^4=\PP H^0(\EEE )$
of curves $(s)_0$ of zeros of all the sections $s$ of $\EEE$
contains at most three distinct curves $C_i\in {\mathcal C}^0_5[2]_p$
and their images $\ell_i$ under $\Phi_p$ are precisely all
the lines in $Y_5$ meeting $\Gamma '$ at $u$.
\end{corollary}

In Section \ref{map-rho} we constructed
a map $\rho_X :\Gamma\lra M_X$
and proved that $\rho_X (\Gamma )=M_X^0$ is an irreducible
component of $M_X$.

\begin{lemma}
$M_X^0=M_X^1$.
\end{lemma}

\begin{proof}
Let $\EEE\in M_X^0$ be generic. The sections of $\EEE$ embed $X$
into the Grassmannian $G=G(2,5)$ of lines in $\PP^4$
and $\EEE\simeq \QQQ_G|_X$. The sections of $\QQQ_G$
correspond to linear forms on $\PP^4$ and the restriction
to $X$ sends them isomorphically onto the sections of $\EEE$.
Let $q\subset X$ be any smooth conic. Considered as
a family of lines in $\PP^4$, $q$ is either 
a pencil of generators of a quadric
in $\PP^3\subset\PP^4$, or the curve of tangents to a conic in $\PP^2$.
In the first case, there is a unique, up to proportionality,
section of $\QQQ_G$ vanishing on $q$, the one defined by the $\PP^3$
which is the linear span of the quadric. In the second case there is
a pencil of such sections in $\PP H^0(\QQQ_G)$. The span of the lines
parametrized by a reducible conic $q$ is $\PP^3$, hence the
same is true for a generic $q$. So, for generic $q$, 
the projective space $\PP H^0(\EEE )$ of zero loci of sections
of $\EEE$ contains a unique reducible curve $D$ having $q$ as one of its
components, and $D$ does not contain other conics or lines.
Hence $D$ is a union of the conic $q$ and of a rational twisted
cubic. Hence $\EEE\in
M_X(q)=M_X^1$.
\end{proof}

$M_X^0$ being birational to $\Gamma$, we obtain:

\begin{corollary}
$\Gamma\simeq\Gamma '\simeq\Gamma ''$.
\end{corollary}

\section{Smoothness and irreducibility of $M_X$}\label{proof-MT}

Let $X=X_{12}$ be a generic Fano threefold of degree 12
with Picard number 1, and $\Gamma$ its dual curve of genus 7.
Let $M_X$, $\rho$ be defined as in Section \ref{map-rho} and let
$M_X^0$ be the component of $M_X$ that is the image of $\rho$.
Propositions \ref{smoothness} and \ref{irredMX} prove Theorem \ref{main}.

\begin{proposition}\label{smoothness}
For generic $X$, $M_X^0$ is in the smooth locus of $M_X$.
\end{proposition}

\begin{proof}
Let $\EEE\in M_X^0$. By Lemma \ref{ncx}, to prove the smoothness
of $M_X$ at $\EEE$ it suffices to find a quintic $C$ in
$\PP H^0(X,\EEE )$ such that $\NNN_{C/X}\not\simeq \OOO_C\oplus\OOO_C(1)$.
Lemmas \ref{movable-bisecants}, \ref{ncx-nodal} and \ref{badline} prove,
for generic $p\in X$, the
existence of such a quintic having a node at $p$.
\end{proof}

\begin{lemma}\label{movable-bisecants}
Let $\EEE\in M_X^0$, $p\in X$ a generic point.
Let $Z_p=\{ C_i\}_{1\leq i\leq n}$ ($n=1,2$ or $3$) be 
the set of curves with singularity at $p$
which are zero loci of sections of $\EEE$, $q_1,\ldots ,q_{24}$
all the conics on $X$ passing through $p$ and $\tilde{C}_i,q'_j$
the proper transforms of the above curves on the blowup
of $p$. Then $\tilde{C}_i\cap q'_j=\empty$ for all $i=1,\ldots ,n$,
$j=1,\ldots ,24$.
\end{lemma}

\begin{proof}
The curves $C_i$ are rational
quintics with only one node. The assertion of Lemma 
is not satisfied for $\EEE ,p$ if and only if there is a pair $i,j$,
for which $C_i$ meets $q_j$ in such a way
that $\length (C_i\cap q_j)\geq 3$.

We want to prove that for generic $p\in X$, there is no conic $q$
through $p$ meeting any one of the curves $C_i$ in such a way
that $\length (C_i\cap q)\geq 3$. Assume the contrary.
By Theorem \ref{bir1}, the number of conics through a generic
point $p$ is 24. Hence, if we assume
that at most seven conics $q_j$ can meet a curve $C_i$ with length $\geq 3$,
then such conics form a proper irreducible component
of $\FFF (X)$. This contradicts the irreducibility of the family
of conics on $X$ (Lemma \ref{conics}).
But it is immediate to see that two different conics $q_1,q_2$
passing through $p$ cannot meet
$C_i$ with length of intersection $\geq 3$. Indeed, if we assume that
such conics $q_1,q_2$ exist, then $C_i\cup q_j$ ($j=1,2$) are
two distinct reducible fibers of the ruled surface $M$
(notations from Theorem \ref{bir1}). The flop
transforms the components $q_j$ into secant lines $l_j$
to the curve $\Gamma\subset Y_5$, and the proper transform
of $C_i$ is contracted to a point of $\Gamma$. This is
absurd, because two distinct fibers of $M$ are contracted
to two distinct points of $\Gamma$. Hence each $C_i$
meets at most one conic $q_j$ and we are done.
\end{proof}

\begin{lemma}\label{ncx-nodal}
Let $\EEE\in M_X^0$, $U\subset X$ the open set of points $p\in X$
for which the birational map $\Phi_p$ of Theorem \ref{bir1}
exists, and $p\in U$. Assume that there are at least two nodal rational
quintic curves $C_1,C_2$ with only one node at $p$ as singularity
which are zero loci of sections of $\EEE$.
Assume also that the proper transform $\tilde{C_i}$ of $C_i$ ($i=1,2$)
on the blowup $X'$ of $X$ at $p$ does not meet any
of the flopping curves of $\kappa$
introduced in Theorem \ref{bir1}. Then the normal
bundle of $C_i$ in $X$ is indecomposable and $M_X$
is nonsingular at the point representing the vector
bundle $\EEE$.
\end{lemma}

\begin{proof}
Assume the contrary, that is, $\NNN_{C/X}\simeq \OOO_C\oplus\OOO_C(1)$,
where $C$ is one of the two curves $C_i$, say, $C=C_1$.
Then $\tilde{C}\simeq\PP^1$ and $\nu^*\NNN_{C/X}\simeq
\OOO_{\PP^1}\oplus\OOO_{\PP^1}(5)$, where $\nu :\tilde{C}\lra C$
is the normalization map; $\nu =\sigma |_{\tilde{C}}$.
Let $\NNN_{\tilde{C}}$ denote the normal bundle to $C$
as a parametrized curve in $X$:\ \ 
$\NNN_{\tilde{C}}=\coker (d\nu :\TTT_{\tilde{C}}\lra \nu^*\NNN_{C/X})$.
Then $\nu^*\NNN_{C/X}$ is obtained from $\NNN_{\tilde{C}}$
by a positive elementary transformation described in \cite{HH}:
$\nu^*\NNN_{C/X}=\elm^+_T\NNN_{\tilde{C}}$, where $T$
is the set of two points in $\PP (\NNN_{\tilde{C}})$
corresponding to the tangent directions of the two branches
at $p$. Hence $\NNN_{\tilde{C}}\simeq
\OOO_{\PP^1}(-2)\oplus\OOO_{\PP^1}(5)$,
$\OOO_{\PP^1}(-1)\oplus\OOO_{\PP^1}(4)$ or
$\OOO_{\PP^1}\oplus\OOO_{\PP^1}(3)$. Hence
$\NNN_{\tilde{C}/X'}\simeq \NNN_{\tilde{C}}(-T)\simeq 
\OOO_{\PP^1}(-4)\oplus\OOO_{\PP^1}(3)$,
$\OOO_{\PP^1}(-3)\oplus\OOO_{\PP^1}(2)$ or
$\OOO_{\PP^1}(-2)\oplus\OOO_{\PP^1}(1)$. 

In the notations of
Theorem \ref{bir1}, the image $\ell$
of $\tilde{C}$ under $\tau\circ\kappa$ is a line in $Y$
meeting $\Gamma$ transversely at one point $z=z(p,\EEE )$, so we have for
the normal
bundles $\NNN_{\ell /Y}=\elm^+_{\mathrm{1\ point}}\NNN_{\tilde{C}/X'}$
Hence the only possible case is
$\NNN_{\ell /Y}\simeq \OOO_{\PP^1}(-1)\oplus\OOO_{\PP^1}(1)$.
According to the description of lines on $Y$,
if $z\not\in C_6^0$ (notations from the proof of Lemma \ref{quintp}),
then there is another line on $Y$ passing through $z$ with normal bundle
$\OOO_{\PP^1}\oplus\OOO_{\PP^1}$ giving the other nodal
quintic $C_2$ and, in replacing $C$ by $C_2$ in the above
argument, we come to the conclusion that
$\NNN_{C_2/X}\not\simeq \OOO_{C_2}\oplus\OOO_{C_2}(1)$. This ends the proof.
\end{proof}

We have seen that $\EEE$
may be a singular point of $M_X$ only if for generic $p\in X$
$C$ is the unique nodal quintic from $\PP H^0(X,\EEE )$, or, equivalently,
$z(p,\EEE )\in C_6^0$ and $\ell =T_{z}C_6^0$ is the unique line on
$Y$ passing through $z(p,\EEE )$.

\begin{lemma}\label{badline}
Let $\EEE\in M_X^0$, $p\in X$ a generic point.
Then there are at least two different curves
with a singular point at $p$
which are zero loci of sections of $\EEE$.
\end{lemma}

\begin{proof}
The assertion that we are to prove now is equivalent to the fact that
for generic $X$ and generic $p\in X$ the curve
$\Gamma\subset Y=Y_5$ defined by Theorem \ref{bir1}
does not meet the closed orbit $C_6^0$ of $SL(2)$.
We prove this fact in three steps in Lemmas \ref{DP1}, \ref{DP2}
and Proposition \ref{inv-bir1}: we prove, first,
that the family of canonical
curves $\Gamma$ in $Y$ is irreducible. 
Second, we present a canonical curve $\Gamma\in Y$
which does not meet $C_6^0$, and this implies that the generic
$\Gamma$ does not meet $C_6^0$ either. Third,
we show that a generic canonical curve $\Gamma\in Y$ is obtained
as in Theorem \ref{bir1} via the map $\Phi_p$
applied to some $X$ and some $p\in X$.
\end{proof}

\begin{lemma}\label{DP1}
Let $Y$ be a Del Pezzo threefold of degree $5$, that is,
a nonsingular threefold linear section $\PP^6\cap G(2,5)$
of the Grassmannian $G(2,5)$ in $\PP^9$. Then the family
of canonical curves $\Gamma$ of genus $7$ in $Y$ is irreducible
of dimension $24$.
\end{lemma}

\begin{proof}
When saying that $\Gamma$ is {\em a canonical curve in} $Y$, we mean that 
the restriction of the class of the hyperplane section of $Y$ to $\Gamma$
is canonical and that the dimension of the linear span of $\Gamma$ in $\PP^9$
coincides with that of the canonical linear system of $\Gamma$,
so $\langle \Gamma \rangle=\langle Y \rangle=\PP^6$.
By \cite{I1}, $Y$ is unique up to projective equivalence
and $h^0(\PP^6,\III_Y(2))=5$.
Hence the hypotheses of Lemmas \ref{curve-case}, \ref{petri-inj}
are verified for the embedding $\Gamma\into Y\into G(2,5)$
and this embedding is determined uniquely, up to the
projective equivalence, by a point $p\in X$, where $X$
is the dual Fano threefold of $\Gamma$. Let $I$ be
the variety parametrizing all the canonical embeddings $\Gamma\into G(2,5)$
such that $\Gamma$ is a generic genus-7 curve and
$\langle \Gamma \rangle\cap  G(2,5)$ is a (smooth) Del Pezzo threefold.
Then $I$ is irreducible of dimension
$$
\dim I = \dim \mathfrak M_7+\dim X+ \dim PGL(5)= 18+3+24=45
$$
Let $\pr$ be the natural projection from $I$ to the open subset $\UUU$
of $G(7,10)$
parametrizing linear sections of $G(2,6)$ which are (smooth) Del Pezzo
threefolds. Then for generic $u\in\UUU$, the fiber
$\pr^{-1}(u)$ is equidimensional and $\dim \pr^{-1}(u)=\dim I-\dim G(7,10)=24$.
By the monodromy argument, to prove the assertion of the lemma,
it suffices to present a distinguished component of the
family of canonical curves in $Y$. Then it has to be the unique one.

There is such an obvious component: the one containing all the
curves $\Gamma$ obtained from the generic pairs
$(p,X)$ via the map $\Phi_p$ of Theorem \ref{bir1}.
This ends the proof of the lemma.
\end{proof}

\begin{lemma}\label{DP2}
Let $F\subset Y$ be any curve. Then there exists
a reducible smoothable canonical curve
$C_0=E_1\cup E_2\cup q$
such that the following conditions are verified:

(a) $E_1,E_2$ are elliptic quintics meeting each other
transversely at $3$ distinct points;

(b) $q$ is a conic meeting each one of the curves
$E_i$ transversely at $2$ points;

(c) $C_0\cap F=\empty$.

\end{lemma}

\begin{proof}
Choose a conic $q$ not meeting $F$; let $\PP^2_q=\langle q\rangle$.
Choose a generic $3$-secant $\PP_2$ to $Y$ meeting $\PP^2_q$. Let
$z=\PP^2\cap \PP^2_q$, $\{ p_1,p_2,p_3\} =\PP^2_q\cap Y$.
Let $l_1,l_2$ be two generic lines in $\PP^2_q$ passing through
$z$, and $\{ p_4^i,p_5^i\}=l_i\cap q$. Let $\PP^3_i= 
\langle \PP^2,l_i\rangle$ and $\PP^4_i$ generic 4-spaces
containing $\PP^3_i$ ($i=1,2$). Then $E_i=\PP^4_i\cap Y$
are elliptic curves such that $\PP^3_i\cap E_i=
\{ p_1,p_2,p_3,p_4^i,p_5^i\}$, that is, the intersection of
$E_i$ with the residual curve $E_{3-i}\cup q$ defines the
divisor from the linear system $\OOO (1)$. For $q$,
the degree of the intersection with $E_1\cup E_2$ is $4$,
which is $\deg K_q+\deg \OOO (1)|_q$. Hence
$C_0=E_1\cup E_2\cup q$ is embedded by a subsystem of the canonical
system, and since $\langle C_0 \rangle =\PP^6$, it is a canonical
curve. It is obvious that for a generic choice of
the above $\PP^4_i$, we have $\PP^4_i\cap F=\empty$.

It remains to see that $C_0$ is smoothable, but this follows
immediately by the technique of \cite{HH} and from
the known normal bundles to the components of our
curve: $\NNN_{E_i}=2\OOO (1)|_{E_i}$ and $\NNN_q=2\OOO_{\PP^1} (1)$.
\end{proof}

\begin{proposition}\label{DP3}\label{inv-bir1}
Let $Y$ be a Del Pezzo threefold of degree $5$
and $\Gamma$ a generic canonical curve of genus $7$ in $Y$. 
Let $X$ be the dual Fano threefold of $\Gamma$. Then there exists
a point $p\in X$ such that the embedding $\Gamma\into Y$
is obtained via the birational map $\Phi_p$ of Theorem \ref{bir1}.
\end{proposition}

\begin{proof}
We have to invert the construction of the map $\Phi_p$
for generic $\Gamma\subset Y$.
Let $\Gamma_0\into Y$ be an embedding corresponding to some
pair $(p_0,X_0)$. Then the anticanonical
linear system $|-K_{\tilde{Y_0}}|$ on the blowup
$\tilde{Y_0}$ of $\Gamma_0$ in $Y$ defines a small contraction $\pi_0$ of
24 \ $(-1,-1)$-lines which are bisecant to $\Gamma_0$ whose image is
a quartic threefold in $\PP^4$ with 24 ordinary double points, 
and the flop over this small contraction gives rise to
an exceptional $\PP^2$ with normal bundle $\OOO(-1)$
which can be contracted to a nonsingular point.
It is easy to verify that this situation is stable under small
deformations. Let $\Gamma_t\into Y$ be a deformation of
$\Gamma_0\into Y$ with parameter $t$.  By the semicontinuity
of $h^i(\tilde{Y_t}, \OOO (-K_{\tilde{Y_t}}))=0$ and because the
$h^i$ vanish for $t=0$ and $i>0$, the linear system
$|-K_{\tilde{Y_t}}|$ defines also a map $\pi_t$ to $\PP^4$
for small $t$. The base-point-free condition is open,
so $\pi_t$ is a morphism. By the semicontinuity
of the degree of a morphism, $\pi_t$ is birational.
The stability of the $(-1,-1)$-curves under deformations
is well-known. See for example \cite{Kod}, where the stability
of nonsingular subvarieties with negative normal bundle is proved,
and the type of the normal bundle is preserved in our case,
because a vector bundle of type $\OOO (k)\oplus\OOO (k)$
has no nontrivial infinitesimal deformations on $\PP^1$.
So, $\tilde{Y_t}$ contains 24 $(-1,-1)$-curves, which have to be contracted
by the anticanonical system. Their intersection indices with the exceptional
divisor of $\tilde{Y_t}\lra Y$ being constant, they descend
to bisecant lines of $\Gamma_t$ in $Y$.

Now we can blow up all the 24 $(-1,-1)$-curves
and blow down the obtained quadrics $\PP^1\times\PP^1$
along the second ruling simultaneously in all
the varieties $\tilde{Y_t}$ for small $t$.
We will obtain the family of nonsingular
projective threefolds $X'_t$. The fiber $X_0$
possesses an exceptional
$\PP^2$. Again by \cite{Kod}, this $\PP^2$ deforms
in a unique way to a compact submanifold in the neighboring
fibers; it is again $\PP^2$ with normal bundle $\OOO(-1)$, because both
$\PP^2$ and the line bundle $\OOO (k)$ on it have no nontrivial
infinitesimal deformations. Thus there is a relative
contraction $X'_t\lra X_t$ blowing down  $\PP^2$ to a nonsingular point
in each fiber.
The fiber $X_0$ is a Fano
threefold of genus 7 with Picard group $Z$, and these
properties are stable under small deformations.
\end{proof}

\begin{proposition}\label{irredMX}
For a generic Fano threefold $X$ of degree $12$ with Picard
group $\ZZ$, the moduli space
$M_X$ is irreducible.
\end{proposition}

\begin{proof}
Let $\EEE\in M_X$. By Proposition \ref{MXglobgen},
the sections of $\EEE$ define a regular map from
$j:X\lra G(2,5)$. Restricting to a generic hyperplane
section $H$ and using Proposition \ref{irredK3},
we see that $i=j|_H$ embeds $H$ linearly into $G(2,5)$ and
satisfies the hypotheses of Lemma \ref{F8-1}.
Hence the linear subspace $U$ defined in the lemma
is a maximal isotropic subspace in 
$V=H^0(\PP^7,\III_H(2))=H^0(\PP^8,\III_X(2))$.
Hence, by Lemma \ref{F8}, $j$ is a linear embedding associated to
a point $w\in\check{X}$, and hence $\EEE\in M_X^0$.
\end{proof}

\vspace{0.5cm}

\end{document}